\documentclass[12pt]{article}

\usepackage{arydshln}
\usepackage{mathrsfs}
\usepackage{amsmath}
\usepackage{amsfonts}
\usepackage{amssymb}
\usepackage[all,cmtip]{xy}
\usepackage[colorlinks]{hyperref}
\usepackage{color}
\usepackage{cleveref}
\usepackage{ulem}
\usepackage[abs]{overpic}
\usepackage{float}
\usepackage{enumitem}
\usepackage{tikz}

\usepackage{etoolbox}
\patchcmd{\thebibliography}
  {\settowidth}
  {\setlength{\itemsep}{0pt plus 0.1pt}\settowidth}
  {}{}
\apptocmd{\thebibliography}
  {\small}
  {}{}

\usepackage[all,cmtip]{xy}
\usepackage{tikz-cd}
\usepackage[toc,page]{appendix}
\usepackage{graphicx}
\usepackage{subcaption}
\usepackage[utf8]{inputenc}
\usepackage[english]{babel}
\usepackage{dsfont}
\usepackage{xcolor}
\usepackage{arydshln}
\usepackage{faktor}
\usepackage{mathtools}
\usepackage{scalerel}
\usepackage{booktabs}

\def\endpf{\hbox{\vrule height1.5ex width.5em}}

\usepackage{dsfont}
\newcommand\xdashmapsto[2][]{\mathrel{\mapstochar\xdashrightarrow[#1]{#2}}}

\makeatletter
\newcommand*{\da@rightarrow}{\mathchar"0\hexnumber@\symAMSa 4B }
\newcommand*{\da@leftarrow}{\mathchar"0\hexnumber@\symAMSa 4C }
\newcommand*{\xdashrightarrow}[2][]{%
	\mathrel{%
		\mathpalette{\da@xarrow{#1}{#2}{}\da@rightarrow{\,}{}}{}%
	}%
}
\newcommand{\xdashleftarrow}[2][]{%
	\mathrel{%
		\mathpalette{\da@xarrow{#1}{#2}\da@leftarrow{}{}{\,}}{}%
	}%
}
\newcommand{\xdashdownarrow}[2][]{%
	\mathrel{%
		\mathpalette{\da@xarrow{#1}{#2}\da@downarrow{}{}{\,}}{}%
	}%
}
\newcommand*{\da@xarrow}[7]{%
	\sbox0{$\ifx#7\scriptstyle\scriptscriptstyle\else\scriptstyle\fi#5#1#6\m@th$}%
	\sbox2{$\ifx#7\scriptstyle\scriptscriptstyle\else\scriptstyle\fi#5#2#6\m@th$}%
	\sbox4{$#7\dabar@\m@th$}%
	\dimen@=\wd0 %
	\ifdim\wd2 >\dimen@
	\dimen@=\wd2 %
	\fi
	\count@=2 %
	\def\da@bars{\dabar@\dabar@}%
	\@whiledim\count@\wd4<\dimen@\do{%
		\advance\count@\@ne
		\expandafter\def\expandafter\da@bars\expandafter{%
			\da@bars
			\dabar@ 
		}%
	}%
	\mathrel{#3}%
	\mathrel{%
		\mathop{\da@bars}\limits
		\ifx\\#1\\%
		\else
		_{\copy0}%
		\fi
		\ifx\\#2\\%
		\else
		^{\copy2}%
		\fi
	}%
	\mathrel{#4}%
}
\makeatother

\usepackage{graphicx}

\newcommand\tikzmark[1]{%
  \tikz[remember picture,overlay]\coordinate (#1);}

\newcommand{\underbracedmatrixll}[2]{%
  \left(\;\hspace{-.27in}
  \smash[b]{\underbrace{
    \begin{matrix}#1\end{matrix}
  }_{#2}}
  \;\right.
  \vphantom{\underbrace{\begin{matrix}#1\end{matrix}}_{#2}}
}

\newcommand{\underbracedmatrixrr}[2]{%
  \left. \;
  \smash[b]{\underbrace{
    \begin{matrix}#1\end{matrix}
  }_{#2}}
  \;\hspace{-.32in}\right)
  \vphantom{\underbrace{\begin{matrix}#1\end{matrix}}_{#2}}
}

\newtheorem{theorem}{Theorem}[section]
\newtheorem{lemma}[theorem]{Lemma}
\newtheorem{corollary}[theorem]{Corollary}
\newtheorem{proposition}[theorem]{Proposition}

\newtheorem{definition}[theorem]{Definition}
\newtheorem{example}[theorem]{Example}
\newtheorem{remark}[theorem]{Remark}

\newcommand{\Addresses}{{
		\bigskip
		\footnotesize
  
         Hanlong Fang, \par\nopagebreak
        \textsc{School of Mathematical Sciences, Peking University, Beijing, 100871, China.}\par\nopagebreak
         \textit{E-mail address:} \href{mailto:hlfang@pku.edu.cn)}{hlfang@pku.edu.cn}
         \medskip
		
	    Xian Wu, \par\nopagebreak        \textsc{School of Mathematical Sciences, Peking University, Beijing, 100871, China.}\par\nopagebreak         \textit{E-mail address:} \href{mailto:xianwu.ag@gmail.com)}{xianwu.ag@gmail.com}         \medskip

}}

\begin{document}
\medskip
\title{\bf Canonical blow-ups of Grassmannians I: How canonical is a Kausz compactification?}

\author {Hanlong Fang\thanks{$^\dagger$Supported by National Key R\&D Program of China under Grant No.2022YFA1006700 and NSFC-12201012.}\,\, and Xian Wu\thanks{Previously also supported by OPUS grant National Science Centre, Poland grant UMO-2018/29/BST1/01290.}}


\date{}

\maketitle
\begin{abstract}
In this paper, we develop a simple uniform picture incorporating the Kausz compactifications and the spaces of complete collineations by blowing up Grassmannians $G(p,n)$ according to a  torus action $\mathbb G_m$. We show that each space of complete collineations is isomorphic to any maximal-dimensional connected component of the $\mathbb G_m$-fixed point scheme of a Kausz-type compactification. We prove that the Kausz-type compactification is the total family over the Hilbert quotient $G(p,n)/ \! \! / \mathbb G_m$ which is isomorphic to the space of complete collineations. In particular, the Kausz compactifications are toroidal embeddings of general linear groups in the sense of Brion-Kumar.
We also show that the Kausz-type compactifications resolve the Landsberg-Manivel birational maps from projective spaces to Grassmannians, by comparing Kausz's construction with ours. As an application, by studying the foliation we derive resolutions of certain birational maps among projective bundles over Grassmannians. The results in this paper are partially taken from the first author's earlier arxiv post (Canonical blow-ups of grassmann manifolds, arxiv:2007.06200), which has been  revised and expanded in collaboration with the second author. 

\end{abstract}

\section{Introduction}
The spaces of complete collineations are classical objects in algebraic geometry compactifying spaces of linear maps of maximal rank. The history dates back at least to Study \cite {Stu}, Severi \cite{Sev2}, Van der Waerden \cite{Van}, etc. on complete conics from the perspective of enumerative geometry, and later to Semple \cite{Se2}, Alguneid \cite {Al}, Tyrrell \cite{Ty}, etc. in higher dimensions, where the studies are through blowing up determinantal varieties and explicit parametrizations based on the Gaussian elimination. Vainsencher \cite{Va} and Laksov \cite{La0} described the spaces of complete collineations as the closures of rational maps 
via the Pl\"ucker coordinates of Grassmannians, and showed that they are wonderful in the sense of De Concini-Procesi \cite{DP}. 
Recall that given a connected reductive algebraic group $G$,  a $G$-variety is called wonderful if it is smooth, projective, its $G$-stable prime divisors $D_1,\cdots,D_r$ are smooth with simple normal crossings, and the strata
$\cap_{i\in I}D_i$ one-to-one correspond to the closures of the $G$-orbits for any $I\subset\{1,\cdots,r\}$. 


Wonderful varieties appear in classical enumerative geometry as an important tool (see \cite{DP2,TK,LLT}, etc.), and play a central role in the theory of spherical varieties (see \cite{Kn,L,BP}).  Among the most important ones the wonderful compactifications of adjoint semisimple groups have applications in arithmetic geometry and representation theory (see \cite{Fa,Sp1,Lu2,He2}, etc.). However, it is well-known that such nice compactification could not exist for reductive groups not of adjoint type. In 2000, Kausz \cite{Ka} made a significant breakthrough by compactifying general linear groups nicely in analogy with the wonderful compactifications. His construction is by blowing up projective spaces iteratively along determinantal varieties both near zero and infinity in a delicate order. 

The aim of the paper is to further investigate the beautiful geometric properties of the Kausz compactifications and make certain natural generalizations. In particular, we shall develop a simple and uniform picture incorporating the Kausz compactifications and the spaces of complete collineations via rational maps, from the perspective of the Bia{\l}ynicki-Birula decomposition \cite{Bi}. We would like to highlight that this
incorporation is primarily influenced by the study of foliations on  Grassmannians induced by the Euler vector fields, which is initiated in complex geometry by Siu-Hwang \cite{Siu, Hw1} in  addressing the nondeformability of compact Hermitian symmetric spaces. A desingularization of such foliations encodes the centers of general linear groups which is elusive in the framework of the wonderful compactification. 
Recently, our idea of blowing up Grassmannians along the closures of the stable/unstable manifolds according to the   Bia{\l}ynicki-Birula  decomposition has been applied by Ding \cite{Di} to give a nice resolution of the Landsberg-Manivel birational maps.  
\smallskip


Let us be more precise. Let $E$ be a free $\mathbb Z$-module of rank $n$. Given an integer $p$ such that
$0<p<n$, consider the Grassmannian $Gr(p,E)$ embedded in the projective space $\mathbb P(\bigwedge^pE)$ via the Plucker embedding. The general linear group ${\rm GL}(E)$ acts linearly on $\bigwedge^pE$ and the induced action on $\mathbb P(\bigwedge^pE)$ stabilizes $Gr(p,E)$.
Choose a non-trivial decomposition $E=E_1\oplus E_2$ and set $s=\dim(E_1)$. Then
we have a decomposition
\begin{equation}\label{cd}
\bigwedge^{p}E:=\bigoplus\limits_{k=0}^p\bigwedge^{k}E_1\otimes\bigwedge^{p-k}E_2.
\end{equation}
The subgroup $G:={\rm GL}(E_1)\times {\rm GL}(E_2)\subset {\rm GL}(E)$ stabilizes 
(\ref{cd}) where each summand is an irreducible $G$-module. We derive a rational $G$-equivariant map
\begin{equation}\label{kspe}
\mathcal K(s,p,E):Gr(p,E)\xdashrightarrow{\,\,\,\,\,\,\,\,\,\,\,\,\,\,\,}\mathbb P(\bigwedge^{p}E)\times\prod\limits_{k=0}^p\mathbb P(\bigwedge^{k}E_1\otimes\bigwedge^{n-k}E_2).  
\end{equation}
Here we make the convention that when $\bigwedge^{k}E_1\otimes\bigwedge^{n-k}E_2=0$, $\mathbb P(\bigwedge^{k}E_1\otimes\bigwedge^{n-k}E_2)$ is a point. Denote by $\mathcal T(s,p,E)$ the closure of the image of the rational map $\mathcal K(s,p,E)$. Denote by $\mathcal M(s,p,E)$ the image of $\mathcal T(s,p,E)$ under the natural projection to 
$\prod\nolimits_{k=0}^p\mathbb P(\bigwedge^{k}E_1\otimes\bigwedge^{n-k}E_2)$. 

\begin{definition}
We call $\mathcal T(s,p,E)$ {\it the Kausz-type compactifications} and $\mathcal M(s,p,E)$ {\it the generalized spaces of complete collineations}. When there is no ambiguity, we write $Gr(p,E)$, $\mathcal T(s,p,E)$, $\mathcal K(s,p,E)$, $\mathcal M(s,p,E)$ as $G(p,n)$, $\mathcal T_{s,p,n}$, $\mathcal K_{s,p,n}$, $\mathcal M_{s,p,n}$, respectively.    
\end{definition}


We state our main results as follows.
\begin{theorem}\label{gwond}
The inverse of (\ref{kspe}) extends to a regular morphism \begin{equation}\label{rspn}
R_{s,p,n}:\mathcal T_{s,p,n}\longrightarrow G(p,n).    
\end{equation} 
$\mathcal T_{s,p,n}$ is smooth and projective over ${\rm Spec}\,\mathbb Z$ with an action of $G$. Set
\begin{equation}\label{rank2}
r=\min\{s,n-s,p,n-p\}.
\end{equation} 
Then the complement of the open
$G$-orbit in $\mathcal T_{s,p,n}$ consists of $2r$ smooth prime divisors with simple normal crossings $D^+_1, 
\cdots,D^+_r,D^-_1$, $\cdots,D^-_r$ such that the following holds.
\begin{enumerate}[label={\rm(\Alph*)}]
\item $D^+_1\cong D^-_1\cong\mathcal M_{s,p,n}$ are smooth and projective over ${\rm Spec}\,\mathbb Z$. 

\item There is a $G$-equivariant flat retraction 
\begin{equation*}
\mathcal P_{s,p,n}:\mathcal T_{s,p,n}\longrightarrow D_1^-\cong\mathcal M_{s,p,n}  
\end{equation*} such that the restriction  $\mathcal P_{s,p,n}|_{D^+_1}:D_1^+\rightarrow D_1^-$ 
is an isomorphism and that for $2\leq i\leq r$,
\begin{equation}\label{KGE}\mathcal P_{s,p,n}(D^-_i)=\mathcal P_{s,p,n}(D^+_{r+2-i})=:\check D_i\,\,    
\end{equation}

\item The closures of $G$-orbits in $\mathcal T_{s,p,n}$ are one-to-one given by
\begin{equation}\label{inrule}
\bigcap\nolimits_{\,\,i\in I^+}D^+_i\mathbin{\scaleobj{1.1}{\bigcap}} \bigcap\nolimits_{\,\,i\in I^-}D^-_i
\end{equation} 
for subsets $I^+,I^-\subset\{1,2,\cdots,r\}$
such that $\min(I^+)+\min(I^-)\geq r+2$ with the convention that $\min(\emptyset)=+\infty$.

\item  $D_1^-$ is wonderful with the $G$-stable divisors $\check D_i$,  $2\leq i\leq r$.
\end{enumerate}
\end{theorem}

\begin{theorem}\label{K=T} When $n=2s=2p$, $\mathcal T_{p,p,2p}$ is isomorphic to the Kausz compactification of the general linear group ${\rm GL}_p$.
\end{theorem}

In general, $\mathcal T_{s,p,n}$ are derived from the Kausz compactifications by parabolic induction as follows. Define a $G$-equivariant rational map
$\pi_1:Gr(p,E)\dashrightarrow Gr(p,E_1)$ by sending each $W\in Gr(p,E)$ with $W\cap E_2=0$ to 
$\pi_1(W)\in Gr(p,E_1)$ such that $W+E_2=\pi_1(W)\oplus E_2$. 
Since $\pi_1$ is also defined by the projection $\bigwedge^pE\rightarrow\bigwedge^pE_1$, it extend to equivariant morphisms $f_1:\mathcal T(s,p,E)\rightarrow Gr(p,E_1)$. Viewing $Gr(p,E_1)$ as the quotient of ${\rm GL}(E_1)$ by the stabilizer of a subspace $W\subset E_1$, we 
see that $f_1$ is a locally fibration with fiber at $W$ being $\mathcal T(p,p,W\oplus E_2)$. Similar statements hold if we exchange the role of $E_1$ and $E_2$. Then we have 
\begin{proposition}\label{tfiber}
\begin{enumerate}[label=\rm(\Alph*)]
\item If $p<n-s$ and $p<s$, $\mathcal T_{s,p,n}$ ({resp.}~$\mathcal M_{s,p,n}$) is a locally trivial fibration over $G(p,s)$ with the fiber $\mathcal T_{p,p,n-s+p}$ ({resp.}~$\mathcal M_{p,p,n-s+p}$), and a locally trivial fibration over $G(p,n-s)$ with the fiber  $\mathcal T_{p,p,s+p}$ ({resp.}~$\mathcal M_{p,p,s+p}$).
    
\item If $n-s<p<s$, $\mathcal T_{s,p,n}$ ({resp.}~$\mathcal M_{s,p,n}$) is a locally trivial fibration over  $G(p,s)$ with the fiber $\mathcal T_{n-s,n-s,n-s+p}$ ({resp.}~$\mathcal M_{n-s,n-s,n-s+p}$), and a locally trivial fibration over $G(s+p-n,s)$ with the fiber $\mathcal T_{n-s,n-s,2n-s-p}$ ({resp.}~$\mathcal M_{n-s,n-s,2n-s-p}$).
    
\item If $p=n-s<s$, $\mathcal T_{s,p,n}$ ({resp.}~$\mathcal M_{s,p,n}$) is a locally trivial fibration over $G(p,s)$ with the fiber $\mathcal T_{p,p,2p}$ ({resp.}~$\mathcal M_{p,p,2p}$).
\end{enumerate}
\end{proposition}

Now we naturally embed the spaces of complete collineations into the Kausz-type compactifications.
\begin{theorem}\label{red2} The space of complete collineations associated with linear maps from a $p$-dimensional vector space to an $(n-p)$-dimensional one is isomorphic to $\mathcal M_{p,p,n}$, and is isomorphic to a $G$-stable divisor of $\mathcal T_{p,p,n}$ which is parabolic inducted from the Kausz compactifiation of ${\rm GL}_p$.
\end{theorem}

Given a reductive group $G$ with the corresponding adjoint semisimple group denoted by $G_{\rm ad}$, we call a $G$-equivariant embedding $X$ toroidal if the quotient map $\pi:G\rightarrow G_{\rm ad}$ extends to a morphism from $X$ to the wonderful compactification of $G_{\rm ad}$ (see \cite[\S 6.2]{BM}). By Theorems \ref{gwond}, \ref{red2}, it is clear that
\begin{corollary}
The Kausz compactifications are toroidal embeddings of general linear groups.    
\end{corollary}

Notice that the $\mathbb G_m$-action on $E=E_1\oplus E_2$ defined by \begin{equation}\label{lt}
t\cdot(e_1,e_2)\mapsto(e_1,t e_2),\ \ 
\,\,(e_1,e_2)\in E_1\oplus E_2 
\end{equation} induces a $\mathbb G_m$-action on $Gr(p,E)$. Moreover, we have the $\mathbb G_m$-equivariant diagram
\vspace{-0.05in}
\begin{equation}\label{KF}
\begin{tikzcd}
&\mathcal T_{s,p,n}\arrow{d}[swap]{\hspace{.03in}\mathcal P_{s,p,n}\,\,\,\,\,\,}\arrow{rr}{\hspace{-0.03in}R_{s,p,n}}\,\,&& G(p,n)\\
&D^-_1\cong\mathcal M_{s,p,n}&& \\
\end{tikzcd}\vspace{-20pt}\,
\end{equation}
with respect to the action (\ref{lt}) on $G(p,n)$ and the trivial action on $\mathcal M_{s,p,n}$. Following \cite{Fu,BS,Kap}, we define the Hilbert quotient $G(p,n)/ \! \! / \mathbb G_m$ as follows. Take any open set $U\subset G(p,n)$ such that the orbit closures of points in $U$ form a flat family of subschemes. There is a morphism from $U$ to  the Hilbert scheme ${\rm Hilb}(G(p,n))$ parametrizing
subschemes of $G(p,n)$. Define  $G(p,n)/ \! \! / \mathbb G_m$ to be the closure of the image of $U$ in ${\rm Hilb}(G(p,n))$. Thaddeus \cite{Th1} proved that $\mathcal M_{s,p,n}$ is isomorphic to  $G(p,n)/ \! \! / \mathbb G_m$ over the complex field. 
We show that 
\begin{theorem}\label{moduli}
Base change (\ref{KF}) to an algebraically closed field $\mathbb K$. $\mathcal M_{s,p,n}\cong G(p,n)/ \! \! /\mathbb G_m$ and $\mathcal P_{s,p,n}:\mathcal T_{s,p,n}\longrightarrow\mathcal M_{s,p,n}$ is the corresponding total family. 
\end{theorem}


Note that the two groups of divisors $Y_i, Z_j$ in Property $3$ of \cite{Ka} share very similar properties to that of $D^-_{i}, D^+_{i}$ in Property (C) of Theorem \ref{gwond}. This is not a mere coincidence but inherently related to the Landsberg-Manivel birational maps 
as follows (see Proposition \ref{g3kausz} for the general statement).    
\begin{corollary}\label{2kausz} Denote by $\mathcal {LM}:\mathbb P^{p^2}\dashrightarrow G(p,2p)$ the Landsberg-Manivel birational maps.  Denote by ${\mathcal {KA}}:\mathcal T_{p,p,2p}\longrightarrow\mathbb P^{p^2}$ the blow-up constructed in \cite{Ka}. The following diagram commutes.
\vspace{-0.05in}
\begin{equation*}
\begin{tikzcd}
&&\,\,\mathcal T_{p,p,2p}\arrow{dl}[swap]{\hspace{-.03in}{\mathcal {KA}}}\arrow{dr}{\hspace{-0.03in}R_{p,p,2p}}\,\,&\\
&\mathbb P^{p^2}\arrow[dashed]{rr}{\mathcal {LM}}&& G(p,2p) \\
\end{tikzcd}\vspace{-20pt}\,.
\end{equation*}
Roughly speaking, by adding $Y_i, Z_j$ to $\mathbb P^{p^2}$  and then subtracting $D^-_i, D^+_i$, we obtain $G(p,2p)$.  
\end{corollary}



Notice that the source and the sink associated to the $\mathbb G_m$-action on $\mathcal T_{s,p,n}$ are isomorphic. We then produce explicit birational transformations among projective bundles over Grassmannians as follows. Denote by $Gr(p,E_1)$ the sub-Grassmannian of $Gr(p,E)$
consisting of $p$-planes in $E_1$, and denote by $\mathbb P(N_1)$ the projectivization of the normal bundle $N_1$ of $Gr(p,E_1)$ in $Gr(p,E)$. When $p\leq n-s$, define $Gr(p,E_2)$ in the same way; otherwise,  
define $Gr(p,E_2)$ to be the sub-Grassmannian of $Gr(p,E)$
consisting of $p$-planes $W$ such that $\dim W\cap E_1=p+s-n$. Denote by $\mathbb P(N_2)$  the projectivization of the normal bundle $N_2$ of $Gr(p,E_2)$ in $Gr(p,E)$. Denote by $R_1:D^-_1\rightarrow\mathbb  P(N_1)$ and $R_2:D^+_1\rightarrow\mathbb P(N_2)$ the corresponding blow-ups induced by  $R_{s,p,n}$ (see \S \ref{basicif} for more precise definitions). 
\begin{corollary}\label{pbun}
We have the following commutative diagram.
\vspace{-0.05in}
\begin{equation*}
\begin{tikzcd}
&&\,\,D^-_1\cong\mathcal M_{s,p,n}\cong D^+_1\arrow{dl}[swap]{\hspace{-.03in}R_1}\arrow{dr}{\hspace{-0.03in}R_2}\,\,&\\
&\mathbb P(N_1)\arrow[dashed]{rr}&& \mathbb P(N_2) \\
\end{tikzcd}\vspace{-20pt}\,.
\end{equation*}  
\end{corollary}

\medskip

We now briefly describe the organization of the paper and the basic ideas for the proof.   
In \S \ref{iterated}, we fix some notations, define group actions on $\mathcal T_{s,p,n}$, and realize $\mathcal T_{s,p,n}$ as iterated blow-ups of $G(p,n)$. In  \S \ref{vander}, we define the Mille Cr\^epes coordinate charts to parametrize $\mathcal T_{s,p,n}$ by adapting the traditional Gaussian elimination method, and hence prove the smoothness of $\mathcal T_{s,p,n}$. In \S \ref{foliation}, we embed $\mathcal M_{s,p,n}$ into $\mathcal T_{s,p,n}$ as a stable divisor and define the flat retraction $\mathcal P_{s,p,n}$ by studying the explicit Bia{\l}ynicki-Birula decomposition on  $\mathcal T_{s,p,n}$, and prove Theorems \ref{red2}, \ref{gwond},  \ref{moduli}. In  \S \ref{kcpt}, we investigate the Landsberg-Manivel birational maps, and
prove Theorem \ref{K=T} and Corollary \ref{2kausz}.  We study in \S \ref{basicif} fibration structures on the Kausz compactifications and prove  Proposition \ref{tfiber} and  Corollary \ref{pbun}.




\section{Preliminaries}\label{iterated}


We first introduce certain notations. For a positive integer $m$, we denote by $\mathbb A^{m}$ the affine scheme
${\rm Spec}\,\mathbb Z\left[x_{1},\cdots,x_{m}\right]$, and by $\mathbb P^{m}$ the projective space
${\rm Proj}\,\mathbb Z\left[x_{0},\cdots,x_{m}\right]$. 

Without loss of generality, we assume that $2p\leq n\leq 2s$ in the paper. 

Define an index set
\begin{equation*}
\mathbb I_{p,n}:=\{(i_1,i_2,\cdots,i_p)\in\mathbb Z^p:1\leq i_p<i_{p-1}<\cdots<i_1\leq n\}.
\end{equation*}
For each $I=(i_1,i_2,\cdots,i_p)\in\mathbb I_{p,n}$, we define {\it Pl\"ucker coordinate functions} $P_{I}$ on $$\mathbb A^{pn}:={\rm Spec}\mathbb Z\left[x_{ij}(1\leq i\leq p,1\leq j\leq n)\right]$$ to be the $p\times p$-subdeterminant of $(x_{ij})$ consisting of the $i_1$-th, $i_2$-th, $\cdots$, $i_p$-th columns.
Define 
\begin{equation*}
\mathcal {G}(p,n):=\left\{\mathfrak p\in\mathbb A^{pn}:P_{I}\notin\mathfrak p\,\,{\rm for\,\,a\,\,certain\,\,}I\in\mathbb I_{p,n}\right\}. 
\end{equation*} 
For any partial permutation $\Delta:=(\delta_1,\delta_2,\cdots,\delta_p)$ of $(1,2,\cdots,n)$, define closed subschemes 
\begin{equation}\label{u12}
U_{\Delta}:={\rm Spec}\left(\mathbb Z\left[\cdots,x_{ij},\cdots\right]/(x_{1\delta_1}-1,x_{1\delta_2},\cdots,x_{1\delta_p},\cdots,x_{p\delta_1},\cdots,x_{p\delta_p}-1)\right) 
\end{equation}
with the  embeddings denoted by 
$e_{\Delta}:U_{\Delta}\xhookrightarrow{\,\,\,\,\,\,\,}\mathcal G(p,n)$.    
Identifying $U_{\Delta}$ with their images under $[\cdots,P_{I}\circ e_{\Delta},\cdots]_{I\in\mathbb I_{p,n}}$, we derive the standard atlas of the Grassmannian $G(p,n)$.  Denote by $\pi:\mathcal G(p,n)\rightarrow G(p,n)$ the natural projection. For each $x\in G(p,n)$, we denote by  $\widetilde{x}$ any element in the preimage $\pi^{-1}(x)\subset\mathcal G(p,n)$.

For convenience, we  write the Pl\"ucker embedding $Gr(p,E)\hookrightarrow\mathbb P(\bigwedge^{p}E)$  as 
\begin{equation*}
e:G(p,n)\xhookrightarrow{\,\,\,\,\,\,\,\,\,\,\,\,\,\,}\mathbb P^{N_{p,n}},
\end{equation*}  
where $\mathbb P^{N_{p,n}}$ is the projective space of dimension $N_{p,n}:=\frac{n!}{(n-p)!p!}-1$ with homogeneous coordinates  
$$[\cdots ,z_I,\cdots]_{\,I\in\mathbb I_{p,n}}.$$


We make the convention that $r$ is always referred to (\ref{rank2}) in this paper.
For $0\leq k\leq r$, define index subsets  
\begin{equation*}
\mathbb I_{s,p,n}^{k}:=\big\{(i_1,\cdots,i_p)\in\mathbb Z^p:{\footnotesize 1\leq i_p<\cdots<i_{k+1}\leq s\,;s+1\leq i_{k}<i_{k-1}<\cdots<i_1\leq n}\}.
\end{equation*}
Consider the following linear subspace of $\mathbb {P}^{N_{p,n}}$,
\begin{equation}\label{subl}
\{[\cdots ,z_I,\cdots]_{I\in\mathbb I_{p,n}}\in\mathbb {P}^{N_{p,n}}:z_I=0\,,\,\,\forall I\notin\mathbb I_{s,p,n}^k \},\,\,\,\, 0\leq k\leq r.
\end{equation}
Since the subspace (\ref{subl}) is isomorphic to $\mathbb {P}^{N^k_{s,p,n}}$
with $N^k_{s,p,n}=|\mathbb I_{s,p,n}^{k}|-1$,  for convenience, we denote (\ref{subl}) by $\mathbb {P}^{N^k_{s,p,n}}$ and its homogeneous coordinates by  $$[\cdots ,z_I,\cdots]_{I\in\mathbb I^k_{s,p,n}}.$$  By dropping the coordinates with indices not in $\mathbb I_{s,p,n}^k$, we have the natural projection  
\begin{equation*}
\begin{split}
F_s^k:\mathbb {P}^{N_{p,n}}&\dashrightarrow\mathbb {P}^{N^k_{s,p,n}} \\
[\cdots ,z_I,\cdots]_{I\in\mathbb I_{p,n}}&\xdashmapsto{}[\cdots, z_I,\cdots]_{I\in\mathbb I^k_{s,p,n}}.
\end{split}
\end{equation*}

Under the above notation, it is clear that (\ref{kspe}) takes the form 
\begin{equation*}
\mathcal K_{s,p,n}:=(e, f_s^0,\cdots,f_s^r):\mathbb A^{p(n-p)}\dashrightarrow\mathbb {P}^{N_{p,n}}\times\mathbb {P}^{N^0_{s,p,n}}\times\cdots\times\mathbb {P}^{N^r_{s,p,n}}, \end{equation*}
where for $0\leq k\leq r$,
\begin{equation}\label{fsk}
f_s^k:=F_s^k\circ e.    
\end{equation}
In homogeneous coordinates, 
\begin{equation}\label{bpp}
\begin{split}
\mathcal K_{s,p,n}(x)=\big( [\cdots ,P_I(\widetilde x),\cdots]_{I\in\mathbb I_{p,n}},[\cdots,P_I(\widetilde x),\cdots]_{I\in\mathbb I^0_{s,p,n}},\cdots,[\cdots,P_I(\widetilde x),\cdots]_{I\in\mathbb I^r_{s,p,n}}\big)
\end{split}.
\end{equation}

Recall that the subgroup ${\rm GL}_s\times {\rm GL}_{n-s}$ 
has a natural action on $G(p,n)$  by  the matrix multiplication.  By checking at generic points, we can derive that there is a unique ${\rm GL}_s\times {\rm GL}_{n-s}$-action $\Delta_{s,p,n}$ on $\mathcal T_{s,p,n}$ such that for each $g\in {\rm GL}_s\times {\rm GL}_{n-s}$ the following diagram commutes.
\vspace{-0.05in}
\begin{equation*}
\begin{tikzcd} &\mathcal T_{s,p,n}  \arrow{r}{\Delta_{s,p,n}(g)}&[2em]\mathcal T_{s,p,n}\arrow{d}{{R_{s,p,n}}} \\ &G(p,n)\arrow[leftarrow]{u}{R_{s,p,n}} \arrow{r}{\delta_{s,p,n}(g)}&G(p,n)\\ \end{tikzcd}.\vspace{-20pt} \end{equation*}
In particular, for the algebraic $\mathbb G_m$-action (\ref{lt}) on $G(p,n)$, there is a unique equivariant $\mathbb G_m$-action on $\mathcal T_{s,p,n}$.

\medskip

In the following, we give an alternative construction of $\mathcal T_{s,p,n}$ as iterated blow-ups of $G(p,n)$. For $0\leq k\leq r$, denote by  $\mathcal S_k$ the ideal sheaf of $\mathcal O_{G(p,n)}$ generated by
\begin{equation}\label{sk}
\{e^*z_{I}:{I\in \mathbb I_{s,p,n}^k}\}\,, 
\end{equation}
which defines a subscheme 
\begin{equation*}
S_k:=\{ x\in G(p,n):\,P_{I}(\widetilde x)=0\,\,\forall{I\in \mathbb I_{s,p,n}^k}\}\,. 
\end{equation*}
For any permutation $\sigma$ of $\{0,1,\cdots,r\}$, let $g_0^{\sigma}:Y^{\sigma}_0\rightarrow G(p,n)$ be the blow-up of $G(p,n)$ along $S_{\sigma(0)}$, and inductively define $g^{\sigma}_{i+1}:Y^{\sigma}_{i+1}\rightarrow Y^{\sigma}_{i}$, $0\leq i\leq r-1$, to be the blow-up of $Y^{\sigma}_i$ along $(g^{\sigma}_0\circ g^{\sigma}_{1}\circ\cdots\circ g^{\sigma}_i)^{-1}(S_{\sigma(i+1)})$. All fit	into the following diagram.
\vspace{-.08in}
\begin{equation}\label{sblow}
\small
\begin{tikzcd}
&Y^{\sigma}_r\ar{r}{g^{\sigma}_r}&Y^{\sigma}_{r-1}\ar{r}{g^{\sigma}_{r-1}}&\cdots\ar{r}{g^{\sigma}_1}&Y^{\sigma}_0\ar{r}{g^{\sigma}_0}&G(p,n) \\
&&(g^{\sigma}_0\circ\cdots\circ g^{\sigma}_{r-1})^{-1}(S_{\sigma(r)})\ar[hook]{u}&\cdots&(g^{\sigma}_0)^{-1}(S_{\sigma(1)})\ar[hook]{u}&S_{\sigma(0)}\ar[hook]{u}\\
\end{tikzcd}.\vspace{-15pt}
\end{equation}
It is clear that $Y_r^{\sigma}$ is independent of the choice of $\sigma$ by the following well-known result.
\begin{lemma}\label{stur}Let $X_0$ be an integral scheme of finite type, and $\mathscr{I}_1$, $\mathscr{I}_{2}$, $\cdots, \mathscr{I}_m$  non-zero ideal sheaves on $X_0$. Inductively define  $f_i:X_i\rightarrow X_{i-1}$ to be the blow-up of $X_{i-1}$ with respect to the inverse image ideal sheaf $(f_1\circ f_2\circ\cdots\circ f_{i-1})^{-1}\mathscr{I}_i\cdot \mathcal O_{X_{i-1}}$,  $1\leq i\leq m$. Then, $X_m$ is isomorphic to the blow-up of $X_0$ with respect to the product of ideal sheaves $\mathscr I_1\mathscr I_2\cdots\mathscr I_m$. \end{lemma} 

Noticing that the elements in (\ref{sk}) are homogeneous of the same degree, we can derive that
\begin{lemma}\label{cpc}
For any permutation $\sigma$, there is an isomorphism $\nu_{\sigma}:\mathcal T_{s,p,n}\rightarrow Y^{\sigma}_r$ such that the following diagram commutes.
\vspace{-0.05in}
\begin{equation*}
\begin{tikzcd}
&\mathcal T_{s,p,n}\arrow[dashed,swap]{rd}{\hspace{-0.03in}\mathcal K^{-1}_{s,p,n}} \arrow{rr}{\nu_{\sigma}}&&Y^{\sigma}_{r}\arrow{dl}{\hspace{-.03in}(g^{\sigma}_0\circ\cdots\circ g^{\sigma}_{r})} \\
&&\,\,G(p,n)\,\,&\\
\end{tikzcd}\vspace{-20pt}\,.
\end{equation*}
In particular, there is a morphism $R_{s,p,n}:\mathcal T_{s,p,n}\rightarrow G(p,n)$ extending $\mathcal K^{-1}_{s,p,n}$.
\end{lemma}

It is clear that $R_{s,p,n}$ is also given by the projection 
of $\mathcal T_{s,p,n}$ to the first factor $\mathbb {P}^{N_{p,n}}$ of the ambient space $\mathbb {P}^{N_{p,n}}\times\mathbb {P}^{N^0_{s,p,n}}\times\cdots\times\mathbb {P}^{N^r_{s,p,n}}$.

\section{ Mille Cr\^epes coordinates} \label{vander}
In this section, we will define a smooth atlas for $\mathcal T_{s,p,n}$ up to the ${\rm GL}_s\times {\rm GL}_{n-s}$ action. 
For convenience, we make the convention that 
for any $x\in G(p,n)$, we denote by $P_I(x)$ the functions $P_I(\widetilde x)$, since they are only used in such a way that the choice of $\widetilde x$ is irrelevant. 

\subsection{Coordinate charts for \texorpdfstring{$R_{s,p,n}^{-1}(U_p)$}{ee} when \texorpdfstring{$p=n-s$}{ee}} \label{vanderp}
Our method is a generalization of the traditional one used in  \cite {Stu,Sev2,Van,Se2,Al,Ty} by iteratively summing rank-$1$ matrices (for this reason we shall call them {\it Mille Cr\^epes}).  
To recall the traditional  method,  we separately address the case when $p=n-s$ in this subsection. 

Define an affine open subscheme of $G(p,n)$ by
\begin{equation}\label{up1}
 U_p:=\left\{\left(
    \begin{array}{c;{2pt/2pt}c}
        Z & I_{p\times p} \\
    \end{array}
\right): Z\,\,{\rm  is\,\, a\,\,} p\times s\,\,{\rm matrix}\,\right\}
\end{equation}
with coordinates $Z=(z_{ij})$.
Define an index set 
\begin{equation}
\mathbb J_p:=\left\{\left.\left(
\begin{matrix}
i_1&i_2&\cdots&i_p\\
j_1&j_2&\cdots&j_p\\
\end{matrix}\right) \right\vert_{}\begin{matrix}
(i_1,i_2,\cdots,i_p)\,\, {\rm is\,\,a\,\,permutation\,\,of\,\,}(1,2,\cdots,\,p)\\(j_1,j_2,\cdots,j_p)\,\, {\rm is\,\,a\,\,partial\,\,permutation\,\,of\,\,}(1,2,\cdots,s)
\end{matrix}\right\}.
\end{equation}
Associate to each $\tau=\left(\begin{matrix} i_1&\cdots&i_p\\
j_1&\cdots&j_p\\
\end{matrix}\right)\in\mathbb J_p$ a subscheme $\mathcal A^{\tau}$ of $\mathbb {P}^{N_{p,n}}\times\mathbb {P}^{N^0_{s,p,n}}\times\cdots\times\mathbb {P}^{N^p_{s,p,n}}$, and an isomorphism $J^{\tau}:\mathbb A^{p(n-p)}\rightarrow \mathcal A^{\tau}$ as follows.

Let $\mathbb {A}^{p(n-p)}:={\rm Spec}\,\mathbb Z[\overrightarrow A,\overrightarrow B^1,\cdots,\overrightarrow B^p]$ where 
\begin{equation}\label{bp}
\begin{split}
&\overrightarrow A:=(a_{i_kj_k})_{1\leq k\leq p}, \,\,\,\overrightarrow B^{k}:=\left(\ (\xi^{(k)}_{i_kj})_{1\leq j\leq s,\,j\neq j_1,j_2,\cdots,j_k},\ (\xi^{(k)}_{ij_k})_{1\leq i\leq p,\,i\neq i_1,i_2,\cdots,i_k}\right)\,\,{\rm for}\,\,1\leq k\leq p.  
\end{split}
\end{equation}
Define a morphism $\Gamma^{\tau}:\mathbb A^{p(n-p)}\longrightarrow U_p\hookrightarrow\mathcal G(p,n)$  by
\begin{equation}\label{ngamma}
\left( \sum\nolimits_{k=1}^p\,\left(v_1^k,\cdots,v_p^k\right)^T\cdot\left(w_1^k,\cdots,w_s^k\right)\cdot\prod\nolimits_{t=1}^{k}a_{i_tj_t}\hspace{-0.13in}\begin{matrix}
  &\hfill\tikzmark{c2}\\
  &\hfill\tikzmark{d2}
  \end{matrix}\,\,\,\,I_{p\times p} \right)\,,
\tikz[remember picture,overlay]   \draw[dashed,dash pattern={on 4pt off 2pt}] ([xshift=0.5\tabcolsep,yshift=7pt]c2.north) -- ([xshift=0.5\tabcolsep,yshift=-2pt]d2.south);
\end{equation} 
where for $1\leq k\leq p$ 
\begin{equation}\label{nxi}
v_t^k=\left\{\begin{array}{ll}    \xi^{(k)}_{tj_k} \,\,\, & t\in\{1,\cdots,p\}\backslash\{i_1,\cdots,i_k\}\\
    0& t\in\{i_1,i_2,\cdots,i_{k-1}\} \\
  
    1 \,\,\, &t=i_k\\
    \end{array}\right.,\ 
     w_t^k=\left\{\begin{array}{ll}
    \xi^{(k)}_{i_kt} \,\,\, & t\in\{1,\cdots,s\}\backslash\{j_1,\cdots,j_k\}\\
    0& t\in\{j_1,j_2,\cdots,j_{k-1}\} \\
    1 &t=j_k\\
    \end{array}\right..
\end{equation}
Define a rational map   $J^{\tau}:\mathbb A^{p(n-p)}\dashrightarrow\mathbb {P}^{N_{p,n}}\times\mathbb {P}^{N^0_{s,p,n}}\times\cdots\times\mathbb {P}^{N^p_{s,p,n}}$ by 
$J^{\tau}:=\mathcal K_{s,p,n}\circ \Gamma^{\tau}$ where $\mathcal K_{s,p,n}$ is given by (\ref{bpp}).

\begin{example}
Let $X=G(3,6)$ and $\tau=\left(
\begin{matrix}
1&2&3\\
1&2&3\\
\end{matrix}\right)
$.  Then
\begin{equation*}
\begin{split}
&\overrightarrow A=(a_{11},a_{22},a_{33}),\,\,\,\overrightarrow  B^1=(\xi^{(1)}_{12},\xi^{(1)}_{13},\xi^{(1)}_{21},\xi^{(1)}_{31}),\,\,\,\overrightarrow  B^2=(\xi^{(2)}_{23},\xi^{(2)}_{32}),\,\,\,\overrightarrow  B^3=\emptyset.\\
\end{split}
\end{equation*}
The map $\Gamma^{\tau}:\mathbb A^{9}\rightarrow U_p$ is given by 
\begin{equation*}
\left( \begin{matrix}
a_{11}&a_{11}\xi^{(1)}_{12}&a_{11}\xi^{(1)}_{13}\\
a_{11}\xi^{(1)}_{21}&a_{11}(\xi^{(1)}_{21}\xi^{(1)}_{12}+a_{22})&a_{11}(\xi^{(1)}_{21}\xi^{(1)}_{13}+a_{22}\xi^{(2)}_{23})\\
a_{11}\xi^{(1)}_{31}&a_{11}(\xi^{(1)}_{31}\xi^{(1)}_{12}+a_{22}\xi^{(2)}_{32})&a_{11}(\xi^{(1)}_{31}\xi^{(1)}_{13}+a_{22}(\xi^{(2)}_{32}\xi^{(2)}_{23}+a_{33})))\\
\end{matrix}\hspace{-0.13in}\begin{matrix}
 &\hfill\tikzmark{c2}\\
 \\
 &\hfill\tikzmark{d2}
 \end{matrix}\,\,\,\,I_{3\times 3} \right)\,.
\tikz[remember picture,overlay]   \draw[dashed,dash pattern={on 4pt off 2pt}] ([xshift=0.5\tabcolsep,yshift=7pt]c2.north) -- ([xshift=0.5\tabcolsep,yshift=-2pt]d2.south);
\end{equation*} 
\end{example}

\begin{lemma}\label{em}
The rational map  $J^{\tau}$ extends to an embedding of $\mathbb A^{p(n-p)}$.
\end{lemma} 
{\noindent\bf Proof of  Lemma \ref{em}.}  
Without loss of generality, assume that 
$\tau= \left(
\begin{matrix}
p&p-1&\cdots&1\\
s&s-1&\cdots&s-p+1\\
\end{matrix}\right)$.

We first show that $J^{\tau}$ has a regular extension over $\mathbb A^{p(n-p)}$. 
It suffices to extend the rational maps $f_s^k\circ \Gamma^{\tau}:\mathbb A^{p(n-p)}\dashrightarrow\mathbb {P}^{N^k_{s,p,n}}$. When $k=p$ it is trivial.  
When $k=p-1$,  $f_s^{p-1}\circ \Gamma^{\tau}$ takes the following form in terms of the homogeneous coordinates for $\mathbb {P}^{N^{p-1}_{s,p,n}}$:
\begin{equation}\label{f1}
\begin{split}
&(f_s^{p-1}\circ \Gamma^{\tau})(\overrightarrow A,\overrightarrow B^1,\cdots,\overrightarrow B^p ))=[\cdots,P_{I}(\Gamma^{\tau}(\overrightarrow A,\overrightarrow B^1,\cdots,\overrightarrow B^p )),\cdots]_{I\in\mathbb I_{s,p,n}^{p-1}}\\
=&\bigg[\,\,a_{ps},\,a_{ps}\xi^{(1)}_{p(s-1)},\,\cdots,\, a_{ps}\xi^{(1)}_{p1},\, -a_{ps}\xi^{(1)}_{(p-1)s}, -a_{ps}\left(\xi^{(1)}_{(p-1)s}\xi^{(1)}_{p(s-1)}+a_{(p-1)(s-1)}\right),\,\\
&-a_{ps}\left(\xi^{(1)}_{(p-1)s}\xi^{(1)}_{p(s-2)}+a_{(p-1)(s-1)}\xi^{(2)}_{(p-1)(s-2)}\right),\cdots,\\
&(-1)^{k-1} a_{ps}\xi^{(1)}_{(p+1-k)s},\,(-1)^{k-1} a_{ps}\left(\xi^{(1)}_{(p+1-k)s}\xi^{(1)}_{p(s-1)}+a_{(p-1)(s-1)}\xi^{(2)}_{(p+1-k)(s-1)}\right),\,\\
&(-1)^{k-1} a_{ps}\left(\xi^{(1)}_{(p+1-k)s}\xi^{(1)}_{p(s-2)}+a_{2(s-1)}\left(\xi^{(2)}_{(p+1-k)(s-1)}\xi^{(2)}_{(p-1)(s-2)}+a_{(p-2)(s-2)}\xi^{(3)}_{(p+1-k)(s-2)}\right)\right),\\
&\cdots,(-1)^{p-1} a_{ps}\xi^{(1)}_{1s},\,(-1)^{p-1} a_{ps}\left(\xi^{(1)}_{1s}\xi^{(1)}_{p(s-1)}+a_{(p-1)(s-1)}\xi^{(2)}_{1(s-1)}\right),\,\cdots\bigg].\\
\end{split}
\end{equation}
Canceling $a_{ps}$, we derive the extension of $f_s^{p-1}\circ \Gamma^{\tau}$.

To derive the extension for  $f_s^k\circ \Gamma^{\tau}$, $0\leq k\leq p-1$, we introduce special indices $I_k,I_k^*,  I_{\mu\nu}^k, I_{\mu\nu}^{k*}\in\mathbb I^k_{s,p,n}$ as follows. For $0\leq k\leq p$, define
\begin{equation}\label{I_k}
    I_k:=(s+k,s+k-1,\cdots,s-p+k+1)\,;
\end{equation}
for $1\leq k\leq p-1$, define
\begin{equation*}
    I^*_{k}:=(s+k+1,s+k-1,s+k-2,\cdots,s-p+k+3,s-p+k+2,s-p+k)\,;
\end{equation*}
for $0\leq k\leq p$, $s-p+k+1\leq \mu\leq s$, and $1\leq \nu\leq s-p+k$ define
\begin{equation*}
    I_{\mu\nu}^k:=(s+k,s+k-1,\cdots,\widehat \mu,\cdots, \nu)\,;
\end{equation*}
for $0\leq k\leq p$, $s+1\leq \mu\leq s+k$, and $s+k+1\leq \nu\leq n$ define
\begin{equation*}
    I^{k*}_{\mu\nu}:=(\nu,s+k,s+k-1,\cdots,\widehat \mu,\cdots, s-p+k+1)\,.
\end{equation*}
Here $\widehat{\mu}$ denotes the omission of $\mu$. A direct computation yields that for $0\leq k\leq p$,
\begin{equation*}
P_{I_k}(\Gamma^{\tau}(\overrightarrow A,\overrightarrow B^1,\cdots,\overrightarrow B^p))=(-1)^{k(p-k)}\cdot \prod_{t=1}^{p-k} a^{p-k+1-t}_{(p+1-t)(s+1-t)}\,,
\end{equation*}
and that for any $I\in\mathbb I^k_{s,p,n}$, $0\leq k\leq p$, there is a polynomial $Q_{I}$ in  $\overrightarrow A,\overrightarrow B^1,\cdots,\overrightarrow B^p$ such that \begin{equation*}
P_{I}(\Gamma^{\tau}(\overrightarrow A,\overrightarrow B^1,\cdots,\overrightarrow B^p))=Q_{I}(\overrightarrow A,\overrightarrow B^1,\cdots,\overrightarrow B^p)\prod_{t=1}^{p-k} a^{p-k+1-t}_{(p+1-t)(s+1-t)}.
\end{equation*}
Hence $J^{\tau}$ has an extension over $\mathbb {A}^{p(n-p)}$. We still denote the extension by $J^{\tau}$. 

We next show that $J^{\tau}$ is an embedding. Computation yields that
for $0\leq k\leq p-1$, $\mu=s-p+k+1$, and $1\leq \nu\leq s-p+k$,
\begin{equation*}\label{hol3}P_{I^k_{\mu\nu}}(\Gamma^{\tau}(\overrightarrow A,\overrightarrow B^1,\cdots,\overrightarrow B^p))=(-1)^{k(p-k)} \xi^{(p-k)}_{(k+1)\nu}\prod_{t=1}^{p-k} a^{p-k+1-t}_{(p+1-t)(s+1-t)}\,,
\end{equation*}
and for $1\leq k\leq p-1$, $s+1\leq \mu\leq s+k$, and $\nu=s+k+1$,
\begin{equation*}
\begin{split}
&P_{I_{\mu\nu}^{k*}}(\Gamma^{\tau}(\overrightarrow A,\overrightarrow B^1,\cdots,\overrightarrow B^p))=(-1)^{k(p-k)+s+k+1-\mu}\xi^{(p-k)}_{(\mu-s)(s-p+k+1)}\prod_{t=1}^{p-k} a^{p-k+1-t}_{(p+1-t)(s+1-t)},\\
&P_{I^*_{k}}(\Gamma^{\tau}(\overrightarrow A,\overrightarrow B^1,\cdots,\overrightarrow B^p))=(-1)^{k(p-k)+1}(a_{k(s-p+k)}+\xi^{(p-k)}_{k(s-p+k+1)}\xi^{(p-k)}_{(k+1)(s-p+k)})\prod_{t=1}^{p-k} a^{p-k+1-t}_{(p+1-t)(s+1-t)}\,.
\end{split}
\end{equation*}
Now we compute $f_s^k\circ\Gamma^{\tau}$ as follows. For $1\leq k\leq p-1$, \begin{equation}\label{f2}
\begin{split}&(f_s^k\circ \Gamma^{\tau})(\overrightarrow A,\overrightarrow B^1,\cdots,\overrightarrow B^p ))=[\cdots,P_{I}\big(\Gamma^{\tau}(\overrightarrow A,\overrightarrow B^1,\cdots,\overrightarrow B^p)\big),\cdots]_{I\in\mathbb I_{s,p,n}^k}\,\,\,\,\,\,\,\,\,\,\,\,\,\,\,\,\,\,\,\,\,\,\,\,\,\,\,\,\,\,\,\,\,\,\,\,\,\,\,\,\,\,\,\,\,\,\,\,\,\,\,\,\,\,\,\,\,\,\,\,\,\,\,\,\,\,\,\,\,\,\,\,\,\,\,\,\,\,\,\,\,\,\,\,\,\,\,\,\,\,\,\,\,\,\,\,\,\,\,\,\,\,\,\,\,\,\,\,\,\,\,\,\,\,\,\,\,\,\,\,\,\,\,\,\,\,\,\,\,\,\,\,\,\,\,\,\,\,\,\,\\
=&\bigg[\,(-1)^{k(p-k)} \prod_{t=1}^{p-k} a^{p-k+1-t}_{(p+1-t)(s+1-t)},\,\cdots,\,(-1)^{k(p-k)} \xi^{(p-k)}_{(k+1)(s-p+k)}\prod_{t=1}^{p-k} a^{p-k+1-t}_{(p+1-t)(s+1-t)},\\
&\cdots,(-1)^{k(p-k)+1} \left(a_{k(s-p+k)}+\xi^{(p-k)}_{k(s-p+k+1)}\cdot\xi^{(p-k)}_{(k+1)(s-p+k)}\right)\prod_{t=1}^{p-k} a^{p-k+1-t}_{(p+1-t)(s+1-t)},\cdots\bigg]\,\,\,\,\,\,\,\,\,\,\,\,\,\,\,\,\,\,\,\,\,\,\\
=&\bigg[\,1,\,\cdots,\, \xi^{(p-k)}_{(k+1)(s-p+k)},\,\xi^{(p-k)}_{(k+1)(s-p+k-1)},\,\cdots, \xi^{(p-k)}_{(k+1)1},\cdots,(-1) \xi^{(p-k)}_{k(s-p+k+1)},\cdots,\,\\&(-1)^{2} \xi^{(p-k)}_{(k-1)(s-p+k+1)},\cdots,(-1)^{k} \xi^{(p-k)}_{1(s-p+k+1)},\cdots,-\left(a_{k(s-p+k)}+\xi^{(p-k)}_{k(s-p+k+1)}\xi^{(p-k)}_{(k+1)(s-p+k)}\right),\cdots\bigg]\,.\\\end{split}\end{equation}
For $k=0$ and $s\geq p+1$,
\begin{equation}\label{f3}
[\cdots,P_{I}(\Gamma^{\tau}(\overrightarrow A,\overrightarrow B^1,\cdots,\overrightarrow B^p)),\cdots]_{I\in\mathbb I_{s,p,n}^0}\\
=[\,1,\,\cdots,\,\xi^{(p)}_{1(s-p)},\,\xi^{(p)}_{1(s-p-1)},\,\cdots,\,\xi^{(p)}_{11},\,\cdots]\,.
\end{equation}
Note that $z_{ps}((e\circ\Gamma^{\tau})(\overrightarrow A,\overrightarrow B^1,\cdots,\overrightarrow B^p))=a_{ps}\,.$
We conclude that $J^{\tau}$  is an embedding. 

The proof of Lemma \ref{em} is complete.
\,\,\,\,$\endpf$
\medskip

Denote by $\mathcal A^{\tau}$ the image of $\mathbb A^{p(n-p)}$ under $J^{\tau}$.  It is clear that  $\left\{\left(\mathcal A^{\tau},(J^{\tau})^{-1}\right)\right\}_{\tau\in\mathbb J_p}$  is a system of coordinate charts of $R^{-1}_{s,p,n}(U_p)$ where  $(J^{\tau})^{-1}:\mathcal A^{\tau}\rightarrow \mathbb A^{p(n-p)}$ is the inverse map.    We call them the {\it  Mille Cr\^epes coordinate charts} of $R^{-1}_{s,p,n}(U_p)$. 
To prove that the  Mille Cr\^epes coordinate charts give an atlas for $R^{-1}_{s,p,n}(U_p)$, it remains to show that
\begin{lemma}\label{coor2} $\bigcup_{\tau\in\mathbb J_p}\mathcal A^{\tau}= R_{s,p,n}^{-1}(U_p)$.
\end{lemma}
{\bf\noindent Proof of Lemma \ref{coor2}.} Define a permutation  $\sigma_p$ by $\sigma_p(k)=p-k$, $0\leq k\leq p$.  Restricting $(\ref{sblow})$ to $U_p$, we have 
\vspace{-.05in}
\begin{equation}\label{blow1}
\footnotesize
\begin{tikzcd}
&Y^p_p\ar{r}{g^{\sigma_p}_p}&Y^{p}_{p-1}\ar{r}{g^{\sigma_p}_{p-1}}&\cdots\ar{r}{g^{\sigma_p}_1}&Y^{p}_0\ar{r}{g^{\sigma_p}_0}&U_p \\
&&(g^{\sigma_p}_0\circ\cdots\circ g^{\sigma_p}_{p-1})^{-1}(S_0\cap U_p)\ar[hook]{u}&\cdots&(g^{\sigma_p}_0)^{-1}(S_{p-1}\cap U_p)\ar[hook]{u}&S_{p}\cap U_p\ar[hook]{u}\\
\end{tikzcd}.\vspace{-20pt}
\end{equation}
It is clear that $Y^p_p\cong R_{s,p,n}^{-1}(U_p)$, and for $0\leq k\leq p$,
\begin{equation*}
   Y^p_{k}:=Y_k^{\sigma}\cap (g^{\sigma_p}_0\circ\cdots\circ g^{\sigma_p}_{k})^{-1}(U_p)\subset  G(p,n)\times\mathbb {P}^{N^p_{s,p,n}}\times\cdots\times\mathbb {P}^{N^{p-k}_{s,p,n}}\,\,\,\,. 
\end{equation*}

We next define an atlas of $Y^p_{k}$. For $0\leq k\leq p$, define index sets 
\begin{equation}
\mathbb J_{p,k}:=\left\{\left. \left(
\begin{matrix}
i_1&i_2&\cdots&i_k\\
j_1&j_2&\cdots&j_k\\
\end{matrix}\right)\right\vert_{}{\begin{matrix}
(i_1,i_2,\cdots,i_k)\,\, {\rm is\,\,a\,\,partial\,\,permutation\,\,of\,\,}(1,2,\cdots,\,p)\\(j_1,j_2,\cdots,j_k)\,\, {\rm is\,\,a\,\,partial\,\,permutation\,\,of\,\,}(1,2,\cdots,s)
\end{matrix}}
\right\}.	    
\end{equation}
For each $\tau=\left(\begin{matrix}
i_1&\cdots&i_k\\
j_1&\cdots&j_k\\
\end{matrix}\right)\in\mathbb J_{p,k}$,  take integers $i^*_1<\cdots<i^*_{p-k}$, $j^*_1<\cdots<j^*_{s-k}$ such that 
$$\left\{i_1,\cdots,i_k,i^*_1,\cdots,i_{p-k}^*\right\}=\left\{1,2,\cdots,p\right\}\,\,\,{\rm and}\,\,\,\left\{j_1,\cdots,j_k,j^*_1,\cdots,j_{s-k}^*\right\}=\left\{1,2,\cdots,s\right\}.$$
Let $\mathbb {A}^{p(n-p)}:={\rm Spec}\,\mathbb Z[\overrightarrow A,\overrightarrow B^1,\cdots,\overrightarrow B^k,\overrightarrow B^*_k]$ where $\overrightarrow A$,
$\overrightarrow B^j$,  $1\leq j\leq k$, are defined by (\ref{bp}), and
\begin{equation*}
\overrightarrow B^{*}_k:=\left(x^{(k+1)}_{i^*_1j^*_1},x^{(k+1)}_{i^*_1j^*_2},\cdots, x^{(k+1)}_{i^*_1j^*_{s-k}},x^{(k+1)}_{i^*_2j^*_1},\cdots,x^{(k+1)}_{i^*_2j^*_{s-k}},\cdots,x^{(k+1)}_{i^*_{p-k}j^*_{1}}\cdots,x^{(k+1)}_{i^*_{p-k}j^*_{s-k}}\right).
\end{equation*}
Define a map $\Gamma_{p,k}^{\tau}:\mathbb A^{p(n-p)}\rightarrow U_p$ by
\begin{equation*}
\left(C^*_k+ \sum\nolimits_{m=1}^k\left(v_1^k,\cdots,v_p^k\right)^T\cdot\left(w_1^k,\cdots,w_s^k\right)\cdot\prod\nolimits_{t=1}^{m}a_{i_tj_t}\hspace{-0.11in}\begin{matrix}
  &\hfill\tikzmark{c2}\\
  &\hfill\tikzmark{d2}
  \end{matrix}\,\,\,\,I_{p\times p} \right)\,.
\tikz[remember picture,overlay]   \draw[dashed,dash pattern={on 4pt off 2pt}] ([xshift=0.5\tabcolsep,yshift=7pt]c2.north) -- ([xshift=0.5\tabcolsep,yshift=-2pt]d2.south);
\end{equation*} 
Here  $v_t^k$, $w_t^k$ are defined by (\ref{nxi}), and $C^*_k$ is defined by
\begin{equation}\label{w7}
C^*_k=\prod_{t=1}^{k} a_{i_tj_t}\cdot\left(\begin{matrix}
&\cdots&0&0&\cdots&0&0&\cdots&0&\cdots\\
&\cdots&0&x^{(k+1)}_{i^*_1j^*_1}&\cdots&x^{(k+1)}_{i^*_1j^*_2}&0&\cdots& x^{(k+1)}_{i^*_1j^*_{s-k}}&\cdots\\
&\cdots&0&0&\cdots&0&0&\cdots&0&\cdots\\
&\cdots&0&x^{(k+1)}_{i^*_2j^*_1}&\cdots&x^{(k+1)}_{i^*_2j^*_2}&0&\cdots& x^{(k+1)}_{i^*_2j^*_{s-k}}&\cdots\\
&\cdots&0&0&\cdots&0&0&\cdots&0&\cdots\\
&\ddots&\vdots&\vdots&\ddots&\vdots&\vdots&\ddots&\vdots&\ddots\\
&\cdots&0&0&\cdots&0&0&\cdots&0&\cdots\\
&\cdots&0&x^{(k+1)}_{i^*_{p-k}j^*_1}&\cdots&x^{(k+1)}_{i^*_{p-k}j^*_2}&0&\cdots& x^{(k+1)}_{i^*_{p-k}j^*_{s-k}}&\cdots\\
&\cdots&0&0&\cdots&0&0&\cdots&0&\cdots\\
\end{matrix}\right).\\
\end{equation}
For $0\leq k\leq p$, define a rational map  $J_{p,k}^{\tau}:\mathbb A^{p(n-p)}\dashrightarrow\mathbb {P}^{N_{p,n}}
\times\mathbb {P}^{N^p_{s,p,n}}\times\cdots\times\mathbb {P}^{N^{p-k}_{s,p,n}}$ by
$J^{\tau}_{p,k}:=\left(e,f_s^{p},f_s^{p-1},\cdots,f_s^{p-k}\right)\circ \Gamma_{p,k}^{\tau}$.
It is clear that $J^{\tau}_{p,p}$ is isomorphic to $J^{\tau}$ after reordering the factors $\mathbb {P}^{N^k_{s,p,n}}$ of the ambient space.

Similarly to Lemma \ref{em}, we can show that the rational map  $J^{\tau}_{p,k}$ extends to an embedding of $\mathbb A^{p(n-p)}$ for each $0\leq k\leq p$. Denote by $\mathcal A_{k}^{\tau}$ the image of $\mathbb A^{p(n-p)}$ in $\mathbb {P}^{N_{p,n}}
\times\mathbb {P}^{N^{p}_{s,p,n}}\times\cdots\times\mathbb {P}^{N^{p-k}_{s,p,n}}$  under  $J_{p,k}^{\tau}$. 
Notice that by (\ref{sk}) each $g_k^{\sigma_p}$ is isomorphic to the blow-up with respect to the ideal generated by $x^{(k+1)}_{i^*_1j^*_1},\cdots,x^{(k+1)}_{i^*_1j^*_{s-k}},x^{(k+1)}_{i^*_2j^*_1},\cdots,x^{(k+1)}_{i^*_2j^*_{s-k}},x^{(k+1)}_{i^*_{p-k}j^*_1},\cdots,x^{(k+1)}_{i^*_{p-k}j^*_{s-k}}.$ Then by induction we can conclude that for $0\leq k\leq p$,     
$\bigcup\nolimits_{\tau\in\mathbb J_{p,k}}\mathcal A_{k}^{\tau}=Y^p_k$. 

The proof of Lemma \ref{coor2} is complete. \,\,\,\,$\endpf$

\subsection{Coordinate charts for \texorpdfstring{$R_{s,p,n}^{-1}(U_l)$}{ee} } \label{vanderl}
We next define the Mille Cr\^epes coordinate charts in general up to group actions. Notice that compared to the traditional method, our parametrization is adapted to the standard atlas for Grassmannians instead of that for projective spaces, and hence there are two blocks conducting Gaussian elimination process.

For $0\leq l\leq r$, define affine open subscheme $U_l\subset G(p,n)$ by
\begin{equation}\label{ul}
U_l:=\{\,\,\,\,\underbracedmatrixll{Z\\Y}{s-p+l\,\,\rm columns}
  \hspace{-.45in}\begin{matrix}
  &\hfill\tikzmark{a}\\
  &\hfill\tikzmark{b}  
  \end{matrix} \,\,\,\,\,
  \begin{matrix}
  0\\
I_{(p-l)\times(p-l)}\\
\end{matrix}\hspace{-.11in}
\begin{matrix}
  &\hfill\tikzmark{c}\\
  &\hfill\tikzmark{d}
  \end{matrix}\hspace{-.11in}\begin{matrix}
  &\hfill\tikzmark{g}\\
  &\hfill\tikzmark{h}
  \end{matrix}\,\,\,\,
\begin{matrix}
I_{l\times l}\\
0\\
\end{matrix}\hspace{-.11in}
\begin{matrix}
  &\hfill\tikzmark{e}\\
  &\hfill\tikzmark{f}\end{matrix}\hspace{-.3in}\underbracedmatrixrr{X\\W}{(n-s-l)\,\,\rm columns}\,\,\,\,\}
  \tikz[remember picture,overlay]   \draw[dashed,dash pattern={on 4pt off 2pt}] ([xshift=0.5\tabcolsep,yshift=7pt]a.north) -- ([xshift=0.5\tabcolsep,yshift=-2pt]b.south);\tikz[remember picture,overlay]   \draw[dashed,dash pattern={on 4pt off 2pt}] ([xshift=0.5\tabcolsep,yshift=7pt]c.north) -- ([xshift=0.5\tabcolsep,yshift=-2pt]d.south);\tikz[remember picture,overlay]   \draw[dashed,dash pattern={on 4pt off 2pt}] ([xshift=0.5\tabcolsep,yshift=7pt]e.north) -- ([xshift=0.5\tabcolsep,yshift=-2pt]f.south);\tikz[remember picture,overlay]   \draw[dashed,dash pattern={on 4pt off 2pt}] ([xshift=0.5\tabcolsep,yshift=7pt]g.north) -- ([xshift=0.5\tabcolsep,yshift=-2pt]h.south);
\end{equation}
with coordinates
\begin{equation}\label{ulx}
\begin{split}
&Z:=(\cdots,z_{ij},\cdots)_{1\leq i\leq l,\,1\leq j\leq s-p-l},\ \ \ \ \ \ X:=(\cdots,x_{ij},\cdots)_{1\leq i\leq l,\,s+l+1\leq j\leq n},\\
&Y:=(\cdots,y_{ij},\cdots)_{l+1\leq i\leq p,\,1\leq j\leq s-p-l},\ \ \ W:=(\cdots,w_{ij},\cdots)_{l+1\leq i\leq p,\,s+l+1\leq j\leq n}.
\end{split}    
\end{equation}
We remark that in the remainder of the paper  $U_l$, $0\leq l\leq r$, is always referred to (\ref{ul}).

For $0\leq l\leq r$, we define index sets 
\begin{equation*}
\mathbb J_l:=\left\{\left(
\begin{matrix}
i_1&\cdots&i_{r-l}&\cdots&i_r\\
j_1&\cdots&j_{r-l}&\cdots&j_r\\
\end{matrix}\right)\rule[-.38in]{0.01in}{.82in}\,\,\footnotesize\begin{matrix}
(i_{1},\cdots,i_{r-l})\,\, {\rm is\,\,a\,\,partial\,\,permutation\,\,of\,\,}(l+1,\cdots,\,p)\\
(i_{r-l+1},\cdots,i_{p})\,\, {\rm is\,\,a\,\,permutation\,\,of\,\,}(1,\cdots,\,l)\\
(j_{1},\cdots,j_{r-l})\,\, {\rm is\,\,a\,\,paritial\,\,permutation\,\,of\,\,}(s+l+1,\cdots,\,n)\\(j_{r-l+1},\cdots,j_{r})\,\, {\rm is\,\,a\,\,paritial\,\,permutation\,\,of\,\,}(1,\cdots,\,s-p+l)
\end{matrix}\right\}.
\end{equation*}
For each $\tau=\left(\begin{matrix}
i_1&\cdots&i_r\\
j_1&\cdots&j_r\\
\end{matrix}\right)\in\mathbb J_l$, set $\mathbb {A}^{p(n-p)}:={\rm Spec}\,\mathbb Z[\overrightarrow A,\widetilde X,\widetilde Y,\overrightarrow B^{1},\cdots,\overrightarrow B^{r}]$, where 
\begin{equation*}
\begin{split}
&\overrightarrow A:=\left((b_{i_kj_k})_{1\leq k\leq r-l},(a_{i_kj_k})_{r-l+1\leq k\leq r}\right),\\
&\widetilde X:=(\cdots,x_{ij},\cdots)_{1\leq i\leq l,\,s+l+1\leq j\leq n},\ \ \widetilde 
Y:=(\cdots,y_{ij},\cdots)_{l+1\leq i\leq p,\,1\leq j\leq s-p-l},       
\end{split}
\end{equation*}
for $1\leq k\leq r-l$,
\begin{equation*}
\overrightarrow B^{k}:=\left(b_{i_kj_k}, \ (\xi^{(k)}_{i_kj})_{s+l+1\leq j\leq n,\,j\neq j_1,j_2,\cdots,j_k},\ (\xi^{(k)}_{ij_k})_{l+1\leq i\leq p,\,i\neq i_1,i_2,\cdots,i_k}\right),
\end{equation*}
and for $r-l+1\leq k\leq r$,
\begin{equation*}
\overrightarrow B^{k}:=\left(a_{i_kj_k}, \ (\xi^{(k)}_{i_kj})_{1\leq j\leq s-p+l,\,j\neq j_{r-l+1},j_{r-l+2},\cdots,j_k},\ (\xi^{(k)}_{ij_k})_{1\leq i\leq l,\,i\neq i_{r-l+1},i_{r-l+2},\cdots,i_k}\right). 
\end{equation*}
Define a map $\Gamma_l^{\tau}:\mathbb A^{p(n-p)}\rightarrow U_l$ by
\begin{equation}\label{ws}
\left(
\begin{matrix}
\sum\limits_{k=r-l+1}^r\left(\prod\limits_{t=r-l+1}^{k}a_{i_{t}j_t}\right)\cdot\Xi_k^T\cdot\Omega_k &0_{l\times(p-l)}&I_{l\times l}&\widetilde X\\ \widetilde Y&I_{(p-l)\times(p-l)}&0_{(p-l)\times l}&\sum\limits_{k=1}^{r-l}\left(\prod\limits_{t=1}^{k}b_{i_tj_t}\right)\cdot\Xi_k^T\cdot\Omega_k\\
\end{matrix}\right)\,.     
\end{equation}
Here, for $1\leq k\leq r-l$, $\Xi_k:=
\left(v_{l+1}^k,\cdots,v_p^k\right)$ with
\begin{equation*}
     v_t^k=\left\{\begin{array}{ll}
    \xi^{(k)}_{tj_k} \,\,\,\,\,\,\,\,\,\,\,\,\,\,\,\,\,\,\,\,\,\,\,\,\,\,\,\,\,\, & t\in\{l+1,l+2,\cdots,p\}\backslash\{i_1,i_2,\cdots,i_k\}\,\,\,\,\,\,\,\,\,\,\,\,\,\,\,\,\, \,\,\,\,\,\,\,\,\,\,\,\,\,\,\,\,\,\\
    0& t\in\{i_1,i_2,\cdots,i_{k-1}\}\,\,\,\,\,\,\, \\
  
    1 \,\,\,\,\,\,\,\,\,\,\,\,\,\,\,\,\,\,\,\,\,\,\,\,\,\,\,\,\,\, &t=i_k\,\,\,\,\,\,\,\,\,\,\,\,\,\,\,\,\, \,\,\,\,\,\,\,\,\,\,\,\,\,\,\,\,\,\\
    \end{array}\right.,
\end{equation*}
and $\Omega_k:=
\left(w_{s+l+1}^k,\cdots,w_n^k\right)$ with
\begin{equation*}
     w_t^k=\left\{\begin{array}{ll}
    \xi^{(k)}_{i_kt} \,\,\,\,\,\,\,\,\,\,\,\,\,\,\,\,\,\,\,\,\,\,\,\,\,\,\,\,\,\, & t\in\{s+l+1,s+l+2,\cdots,n\}\backslash\{j_1,j_2,\cdots,j_k\}\,\,\\
    0& t\in\{j_1,j_2,\cdots,j_{k-1}\}\,\,\,\,\,\,\, \\
  
    1 \,\,\,\,\,\,\,\,\,\,\,\,\,\,\,\,\,\,\,\,\,\,\,\,\,\,\,\,\,\, &t=j_k\,\,\,\,\,\,\,\,\,\\
    \end{array}\right.;
\end{equation*}
for $r-l+1\leq k\leq r$,  $\Xi_k:=
\left(v_{1}^k,\cdots,v_l^k\right)$ with
\begin{equation*}
     v_t^k=\left\{\begin{array}{ll}
    \xi^{(k)}_{tj_k} \,\,\,\,\,\,\,\,\,\,\,\,\,\,\,\,\,\,\,\,\,\,\,\,\,\,\,\,\,\, & t\in\{1,2,\cdots,l\}\backslash\{i_{r-l+1},i_{r-l+2},\cdots,i_k\}\,\,\,\,\,\,\,\,\,\,\,\,\,\,\,\,\, \,\,\,\,\,\,\,\,\,\,\,\,\,\,\,\,\,\\
    0& t\in\{i_{r-l+1},i_{r-l+2},\cdots,i_{k-1}\}\,\,\,\,\,\,\, \\
  
    1 \,\,\,\,\,\,\,\,\,\,\,\,\,\,\,\,\,\,\,\,\,\,\,\,\,\,\,\,\,\, &t=i_k\,\,\,\,\,\,\,\,\,\,\,\,\,\,\,\,\, \,\,\,\,\,\,\,\,\,\,\,\,\,\,\,\,\,\\
    \end{array}\right.,
\end{equation*}
and $\Omega_k:=
\left(w_{1}^k,\cdots,w_{s-p+l}^k\right)$ with
\begin{equation*}
     w_t^k=\left\{\begin{array}{ll}
    \xi^{(k)}_{i_kt} \,\,\,\,\,\,\,\,\,\,\,\,\,\,\,\,\,\,\,\,\,\,\,\,\,\,\,\,\,\, & t\in\{1,2,\cdots,s-p+l\}\backslash\{j_{r-l+1},j_{r-l+2},\cdots,j_k\} \,\,\,\,\,\,\,\,\,\,\,\,\,\,\,\,\,\\
    0& t\in\{j_{r-l+1},j_{r-l+2},\cdots,j_{k-1}\}\,\,\,\,\,\,\, \\
  
    1 \,\,\,\,\,\,\,\,\,\,\,\,\,\,\,\,\,\,\,\,\,\,\,\,\,\,\,\,\,\, &t=j_k\,\,\,\,\,\,\,\,\,\,\,\,\,\\
    \end{array}\right..
\end{equation*}

We define a rational map   $J_l^{\tau}:\mathbb A^{p(n-p)}\dashrightarrow\mathbb {P}^{N_{p,n}}\times\mathbb {P}^{N^0_{s,p,n}}\times\cdots\times\mathbb {P}^{N^r_{s,p,n}}$ by $J_l^{\tau}:=\mathcal K_{s,p,n}\circ \Gamma_l^{\tau}$.
Similarly to Lemma \ref{em}, we can prove that the rational map  $J_l^{\tau}$ is an embedding of $\mathbb A^{p(n-p)}$. Denote by $\mathcal A^{\tau}$ the image of $\mathbb A^{p(n-p)}$ under $J_l^{\tau}$.  It is clear that  $\left\{\left(\mathcal A^{\tau},(J_l^{\tau})^{-1}\right)\right\}_{\tau\in\mathbb J_p}$  is a system of coordinate charts of $R^{-1}_{s,p,n}(U_l)$, where $(J_l^{\tau})^{-1}:\mathcal A^{\tau}\rightarrow \mathbb A^{p(n-p)}$ is the inverse map. Similarly to Lemma \ref{coor2}, we can further prove that $\bigcup_{\tau\in\mathbb J_l}\mathcal A^{\tau}= R_{s,p,n}^{-1}(U_l)$, and hence $\left\{\left(\mathcal A^{\tau},(J_l^{\tau})^{-1}\right)\right\}_{\tau\in\mathbb J_l}$  is an atlas for $R^{-1}_{s,p,n}(U_l)$. \begin{definition}
We call $\left\{\left(\mathcal A^{\tau},(J_l^{\tau})^{-1}\right)\right\}_{\tau\in\mathbb J_l}$ the {\it  Mille Cr\^epes coordinate charts} of $R^{-1}_{s,p,n}(U_l)$, or simply the {\it  Mille Cr\^epes coordinates}.   
\end{definition}

\begin{example}
Consider $G(4,8)$ with $s=4$, $l=2$, and $\tau=\left(
\begin{matrix}
3&4&1&2\\
7&8&1&2\\
\end{matrix}\right)\in\mathbb J_2$. Then	\begin{equation*}
\begin{split}
&\overrightarrow  A=(b_{37},b_{48},a_{11},a_{22}),\,\,\,\,\,\,\widetilde X=\left(x_{17},x_{18},x_{27},x_{28}\right),\,\,\,\,\,\,\widetilde Y=\left(y_{31},y_{32},y_{41},y_{42}\right),\\
&\overrightarrow  B^1=(b_{37},\xi^{(1)}_{38},\xi^{(1)}_{47}),\,\,\,\,\,\,\overrightarrow  B^2=\emptyset,\,\,\,\,\,\, B^3=(a_{11},\xi^{(3)}_{12},\xi^{(3)}_{21}),\,\,\,\,\,\,\overrightarrow  B^4=\emptyset.\\
\end{split}
\end{equation*}
The map $\Gamma_2^{\tau}:\mathbb C^{16}\rightarrow U_2$ is given by 
\begin{equation*}
\left(
\begin{matrix}
a_{11}&a_{11}\cdot\xi^{(3)}_{12}&0&0&1&0&x_{17}&x_{18}\\
a_{11}\cdot\xi^{(3)}_{21}&a_{11}\cdot(\xi^{(3)}_{12}\cdot\xi^{(3)}_{21}+a_{22})&0&0&0&1&x_{27}&x_{28}\\
y_{31}&y_{32}&1&0&0&0&b_{37}&b_{37}\cdot\xi^{(1)}_{38}\\
y_{41}&y_{42}&0&1&0&0&b_{37}\cdot\xi^{(1)}_{47}&b_{37}\cdot(\xi^{(1)}_{47}\cdot\xi^{(1)}_{38}+b_{48})\\
	\end{matrix}
	\right).
	\end{equation*}
\end{example}

\begin{proposition}\label{tspnsmooth}
$\mathcal T_{s,p,n}$ is  smooth over ${\rm Spec}\,\mathbb Z$.   
\end{proposition}
{\noindent\bf Proof of  Proposition \ref{tspnsmooth}.}  
Since the singular locus of $\mathcal T_{s,p,n}$ 
is ${\rm GL}_s\times {\rm GL}_{n-s}$-invariant, we conclude Proposition \ref{tspnsmooth} by  the  Mille Cr\^epes coordinate charts constructed above. \,\,\,\,$\endpf$

\section{Explicit Bia{\l}ynicki-Birula decomposition on \texorpdfstring{$\mathcal T_{s,p,n}$}{hh}} \label{foliation}
Bia{\l}ynicki-Birula \cite{Bi} decomposed a smooth projective scheme $X$ under an algebraic torus action into a disjoint union of stable/unstable subschemes fibered over connected components of the fixed point scheme. Exploiting such idea in the Mille Cr\^epes coordinate charts, we investigate the geometry of $\mathcal T_{s,p,n}$ and $\mathcal M_{s,p,n}$ in this section. 

Assume $2p\leq n\leq 2s$ and let $r$ be as in (\ref{rank2}). For $0\leq k\leq r$, define subschemes 
of $G(p,n)$ 
\begin{equation*}
\begin{split}
&\mathcal V_{(p-k,k)} :=\left\{\left. \left(
\begin{matrix}
0&X\\
Y&0\\
\end{matrix}\right)\right\vert_{}\footnotesize\begin{matrix}
X\,\,{\rm is\,\,an\,\,}k\times (n-s)\,\,{\rm matrix\,\,of\,\,rank}\,\,k\,\\
Y\,\,{\rm is\,\,a\,\,}(p-k)\times s\,\,{\rm matrix\,\,of\,\,rank}\,\,(p-k)\\
\end{matrix}
\right\}\,,\\
&\mathcal V_{(p-k,k)}^+:= \left\{\left.\left(
\begin{matrix}
0&X\\
Y&W\\
\end{matrix}\right)\right\vert_{}\footnotesize\begin{matrix}
X\,\,{\rm is\,\,an\,\,}k\times (n-s)\,\,{\rm matrix \,\,of\,\,rank\,\,}k\,\\
Y\,\,{\rm is\,\,a\,\,}(p-k)\times s\,\,{\rm matrix \,\,of\,\,rank\,\,}(p-k)\,\\
\end{matrix}\right\},\\
&\mathcal V_{(p-k,k)}^-:= \left\{\left.\left(
\begin{matrix}
Z&X\\
Y&0\\
\end{matrix}\right)\right\vert_{}{\footnotesize\begin{matrix}
X\,\,{\rm is\,\,an\,\,}k\times (n-s)\,\,{\rm matrix \,\,of\,\,rank\,\,}k\,\\
Y\,\,{\rm is\,\,a\,\,}(p-k)\times s\,\,{\rm matrix \,\,of\,\,rank\,\,}(p-k)\,\\
\end{matrix}}\right\}\,.    
\end{split}
\end{equation*}
A direct computation yields that 
the connected components of the fixed point scheme of $G(p,n)$ under the $\mathbb G_m$-action defined by (\ref{lt}) are
$\mathcal V_{(p,0)},\mathcal V_{(p-1,1)},\cdots,\mathcal V_{(p-r,r)}$. 
Moreover, for $0\leq k\leq r$,
\begin{equation*}
\overline{\mathcal V_{(p-k,k)}^+}\mathbin{\scaleobj{1.2}{\backslash}} \mathcal V_{(p-k,k)}^+=\bigsqcup\nolimits_{j=k+1}^r\mathcal V_{(p-j,j)}^+\,\,\,\,\,\,\,\,{\rm and}\,\,\,\,\,\,\,\overline{\mathcal V_{(p-k,k)}^-}\mathbin{\scaleobj{1.2}{\backslash}} \mathcal V_{(p-k,k)}^-=\bigsqcup\nolimits_{j=0}^{k-1}\mathcal V_{(p-j,j)}^-\,,
\end{equation*}
where $\overline{\mathcal V_{(p-k,k)}^{\pm}}$ is the Zariski closure of $\mathcal V_{(p-k,k)}^{\pm}$ and  $\mathcal V_{(p-r-1,r+1)}^+=\mathcal V_{(p+1,-1)}^-=\emptyset$ by convention. 
\begin{definition}
We call $\mathcal V_{p,0}$ the source and $\mathcal V_{p-r,r}$ the sink of $G(p,n)$ under the $\mathbb G_m$-action (\ref{lt}). 
\end{definition}

We have that
\begin{lemma}\label{kl} Let $S_k\subset G(p,n)$ be the subscheme defined by (\ref{sk}). For $0\leq k\leq r$, $S_k$ is the  scheme-theoretic union of  $\overline{\mathcal V_{(p-k-1,k+1)}^+}$ and  $\overline{\mathcal V_{(p-k+1,k-1)}^-}$. 
\end{lemma}
{\noindent\bf Proof of  Lemma \ref{kl}.} It follows by \cite[Theorem 1.4]{BV}. 
\,\,\,$\endpf$
\smallskip

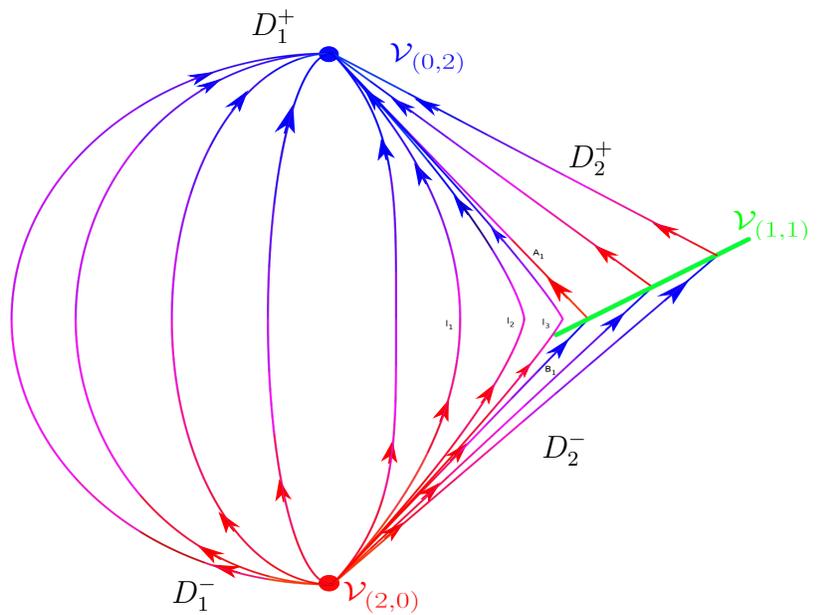
\begin{figure}[htbp] {
\begin{center}
 \begin{overpic}[angle=-90,origin=c,height=11cm,width=16cm]
{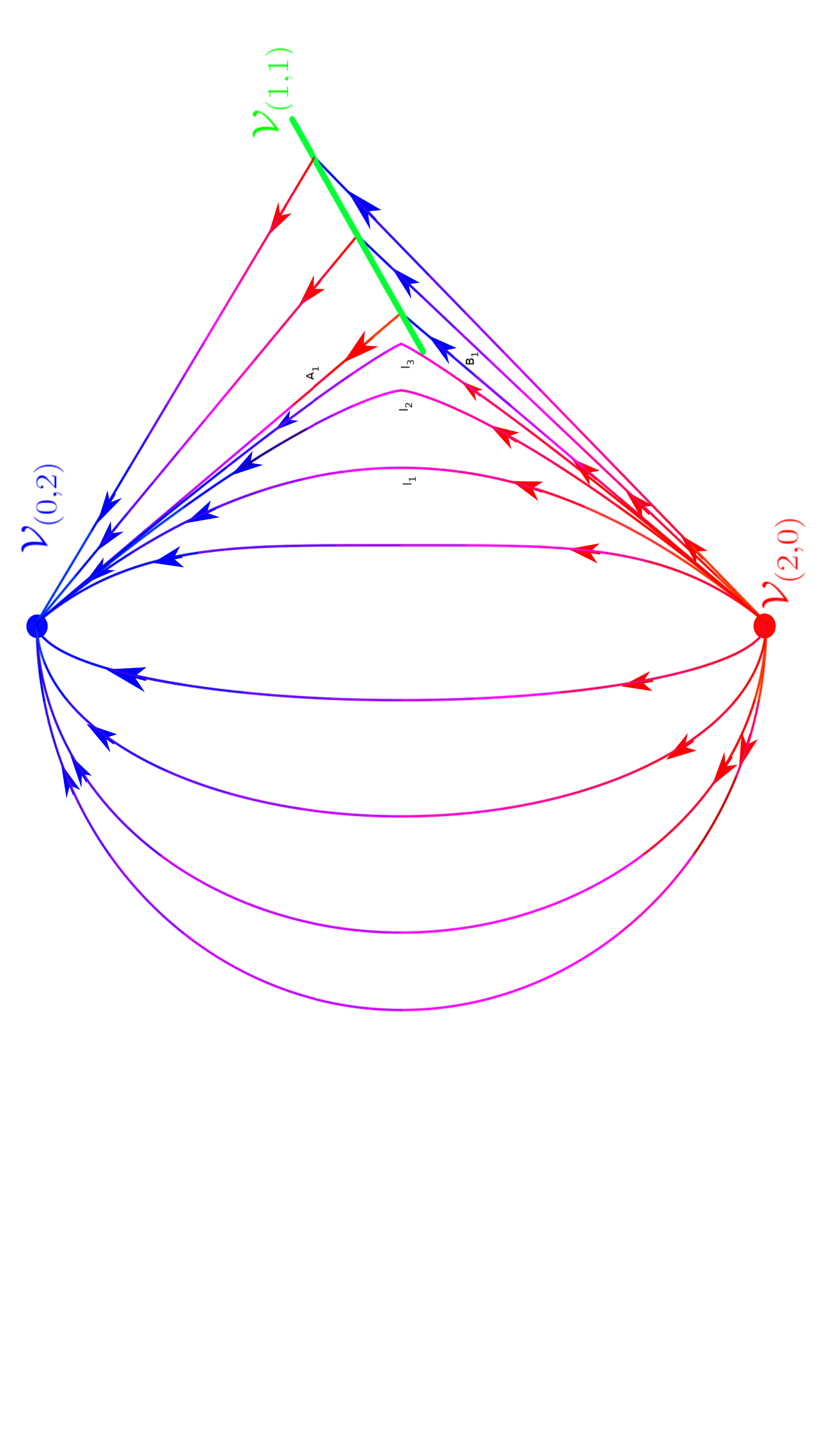}
\put (230,310) {$D^+_1$}
\put (350,260) {$D^+_2$}
\put (200,95) {$D^-_1$}
\put (340,150) {$D^-_2$} 
\end{overpic}
\end{center}}
\vspace{-1.3in}
\caption{The foliation on $G(2,4)$. $D^+_1$ is the exceptional divisor over the infinity point $\mathcal V_{0,2}$. $D^-_1$ is the exceptional divisor over the origin $\mathcal V_{2,0}$. $D^+_2$ is the exceptional divisor over the union of all orbital rational curves between $\mathcal V_{1,1}$ and $\mathcal V_{0,2}$; $D^-_2$ is the exceptional divisor over the union of all orbital rational curves between $\mathcal V_{1,1}$ and $\mathcal V_{2,0}$.}
\label{G2}
\end{figure}

Notice that in each Mille Cr\^epes coordinate chart  
$\mathcal A^{\tau}={\rm Spec}\,\mathbb Z[\widetilde X,\widetilde Y,\overrightarrow B^{1},\cdots,\overrightarrow B^{p}]$ of $R_{s,p,n}^{-1}(U_l)$, the action $\mathbb G_m={\rm Spec}\,\mathbb Z[t,t^{-1}]$ in (\ref{lt}) takes the form $$b_{i_1j_1}\mapsto t b_{i_1j_1},\,\,a_{i_{p-l+1}j_{p-l+1}}\mapsto t^{-1}a_{i_{p-l+1}j_{p-l+1}}.$$ 
Computation yields that the connected components of the fixed point scheme of $\mathcal T_{s,p,n}$ under  $\mathbb G_m$ are given by 
\begin{equation}\label{11}
\mathcal D_{(p-k,k)}:=R_{s,p,n}^{-1}(\mathcal V_{(p-k,k)}), \,\,0\leq k\leq r,
\end{equation}
and the stable/unstable subschemes are given by
\begin{equation}\label{tschubert}
\mathcal D^{\pm}_{(p-k,k)}:=R_{s,p,n}^{-1}(\mathcal V^{\pm}_{(p-k,k)}), \,\,0\leq k\leq r.
\end{equation}
It is clear that $\mathcal D_{(p,0)}$ and $\mathcal D_{(p-r,r)}$ are smooth divisors of $\mathcal T_{s,p,n}$, and  $\mathcal D_{(p-k,k)}$, $1\leq k\leq r-1$,  are codimension two smooth closed subschemes. 
Moreover, for $0\leq k\leq r$,
\begin{equation*}%
\overline{\mathcal D_{(p-k,k)}^+}\mathbin{\scaleobj{1.2}{\backslash}} \mathcal D_{(p-k,k)}^+=\bigsqcup\nolimits_{j=k+1}^r\mathcal D_{(p-j,j)}^+\,\,\,\,\,\,\,\,{\rm and}\,\,\,\,\,\,\,\overline{\mathcal D_{(p-k,k)}^-}\mathbin{\scaleobj{1.2}{\backslash}} \mathcal D_{(p-k,k)}^-=\bigsqcup\nolimits_{j=0}^{k-1}\mathcal D_{(p-j,j)}^-\,.
\end{equation*}
For $1\leq k\leq r$, define
\begin{equation}\label{boundary}
D_{k}^-:=\overline{\mathcal D_{(p-k+1,k-1)}^-},\,\,D_{k}^+:=\overline{\mathcal D_{(p-r+k-1,r-k+1)}^+}.    
\end{equation}
We call $D_1^-$ and $D_1^+$  {\it the source} and {\it the sink} of $\mathcal T_{s,p,n}$ under the $\mathbb G_m$-action induced by (\ref{lt}), respectively. 
See Figure (\ref{G2}) for illustration.

\begin{lemma}\label{snc}
$D_{1}^+$, $\cdots$, $D_{r}^+$, $D_{1}^-$,  $\cdots$, $D_{r}^-$ are distinct smooth divisors with simple normal crossings. In particular,  $D_1^-=\mathcal D_{(p,0)}$,  $D_1^+=\mathcal D_{(p-r,r)}$.
\end{lemma}

{\noindent\bf Proof of  Lemma \ref{snc}.} 
Let $\left(\mathcal A^{\tau},(J_l^{\tau})^{-1}\right)$ be any coordinate chart in the  Mille Cr\^epes coordinates with $0\leq l\leq r$ and $\tau=\left(
	\begin{matrix}
	i_1&i_2&\cdots&i_{r-l}&\cdots&i_r\\
	j_1&j_2&\cdots&j_{r-l}&\cdots&j_r\\
	\end{matrix}\right)\in\mathbb J_l$. Then
\begin{equation*}
\begin{split}
&D_k^+\bigcap \mathcal A^{\tau} =\left\{a_{i_{k}j_{k}}=0\right\}\,\,\,{\rm when}\,\, r-l+1\leq k\leq r\,;D_k^+\bigcap \mathcal A^{\tau} =\emptyset\,\,\,{\rm when}\,\, 1\leq k\leq r-l\,.\\
&D_k^-\bigcap \mathcal A^{\tau} =\left\{b_{i_{k-l}j_{k-l}}=0\right\}\,\,\,{\rm when}\,\, l+1\leq k\leq r\,;D_k^-\bigcap \mathcal A^{\tau} =\emptyset\,,\,\,{\rm when}\,\, 1\leq k\leq l\,.\\
\end{split}
\end{equation*}
Applying the ${\rm GL}_s\times {\rm GL}_{n-s}$-action, we conclude Lemma \ref{snc}. \,\,\,$\endpf$

\begin{remark}\label{nd}
Base change to an algebraically closed field $\mathbb K$.  For any $0\leq k\leq r$ and  closed point $\mathfrak a\in\mathcal D_{(p-k,k)}$, 
decompose $T_{\mathfrak a}(\mathcal T_{s,p,n})=T_{\mathfrak a}(\mathcal T_{s,p,n})^0\oplus T_{\mathfrak a}(\mathcal T_{s,p,n})^{+}\oplus T_{\mathfrak a}(\mathcal T_{s,p,n})^{-}$ into $\mathbb G_m$-weight spaces with respect to the zero, positive, and negative weights. Then 
\begin{enumerate}[label=(\roman*)]
\item when $k=0$, $\dim T_{\mathfrak a}(\mathcal T_{s,p,n})^{+}=1$ and $\dim T_{\mathfrak a}(\mathcal T_{s,p,n})^{-}=0$; 
		
\item when $1\leq k\leq r-1$, $\dim T_{\mathfrak a}(\mathcal T_{s,p,n})^{+}=\dim T_{\mathfrak a}(\mathcal T_{s,p,n})^{-}=1$; 
	
\item when $k=r$, $\dim T_{\mathfrak a}(\mathcal T_{s,p,n})^{+}=0$ and $\dim T_{\mathfrak a}(\mathcal T_{s,p,n})^{-}=1$. 
\end{enumerate}
When $\mathbb K=\mathbb C$, the induced  holomorphic vector field  on $G(p,n)$ takes the form 
\begin{equation*}
-\sum\nolimits_{i=1}^l\sum\nolimits_{j=1}^{s-p+l}z_{ij}\frac{\partial}{\partial z_{ij}}+\sum\nolimits_{i=l+1}^{p}\sum\nolimits_{j=s+l+1}^{n}w_{ij}\frac{\partial}{\partial w_{ij}}\,
\end{equation*}
in a certain coordinate chart around $\mathcal V_{(p-l,l)}$, $0\leq l\leq r$.
In particular, there are points with infinitely many orbital curves passing through. After lifting to $\mathcal T_{s,p,n}$, the corresponding vector field takes the normal form $b\frac{\partial}{\partial b}$ on $\mathcal D_{(p,0)}$, $-a\frac{\partial}{\partial a} +b\frac{\partial}{\partial b}$ on $\mathcal D_{(p,l)}$ for $1\leq l\leq r-1$, and $-a\frac{\partial}{\partial a}$ on $\mathcal D_{(p-r,r)}$. In particular, for each point either one line flows in/out, or one line flows in and one line flows out. This is the original motivation to blow up $G(p,n)$ in such a way. 
\end{remark}



We next construct a flat map   $\mathcal P_{s,p,n}:\mathcal T_{s,p,n}\rightarrow D^-_1$. Let
$\mathfrak {P}:\mathcal T_{s,p,n}\rightarrow \mathbb {P}^{N^0_{s,p,n}}\times\cdots\times\mathbb {P}^{N^r_{s,p,n}}$ be the restriction of the  natural projection 
$\mathbb {P}^{N_{p,n}}\times\mathbb {P}^{N^0_{s,p,n}}\times\cdots\times\mathbb {P}^{N^r_{s,p,n}}\rightarrow\mathbb {P}^{N^0_{s,p,n}}\times\cdots\times\mathbb {P}^{N^r_{s,p,n}}$.

\begin{lemma}\label{emb}
$\mathfrak P$ induces an embedding of $D^-_1$ into $\mathbb {P}^{N^0_{s,p,n}}\times\cdots\times\mathbb {P}^{N^r_{s,p,n}}$.
\end{lemma}
{\bf\noindent Proof of Lemma \ref{emb}.} Notice that $\mathcal V_{(p,0)}$ is isomorphic to a sub-Grassmannian $G(p,s)$ of $G(p,n)$, and the subset of the Pl\"ucker coordinate functions  $\{P_I\}_{I\in\mathbb I^0_{s,p,n}}$ of $G(p,n)$ is the full set of the Pl\"ucker coordinate functions $\{P_I\}_{I\in\mathbb I_{p,s}}$ of $G(p,s)$. Then the rational map  $f^0_s$ given by (\ref{fsk}) is well-defined on $\mathcal V_{(p,0)}$ and $f^0_s\big|_{\mathcal V_{(p,0)}}:\mathcal V_{(p,0)}\rightarrow \mathbb {P}^{N^0_{s,p,n}}$ is an embedding. 

By computing in the Mille Cr\^epes coordinate charts as in the proof of Lemma \ref{em}, we then conclude that the projection of 
$D^-_1$ to $\mathbb {P}^{N^0_{s,p,n}}\times\cdots\times\mathbb {P}^{N^r_{s,p,n}}$ is an embedding. 
\,\,\,\,$\endpf$.

\begin{lemma}\label{iemb}
The image of $\mathcal T_{s,p,n}$ under $\mathfrak P$ is the image of $D^-_1$ under $\mathfrak P$.
\end{lemma}
{\bf\noindent Proof of Lemma \ref{iemb}.} It suffices to show that $\mathfrak P(\mathcal A^{\tau})\subset \mathfrak P\left(D^-_1\right)$, where $\left(\mathcal A^{\tau},(J_0^{\tau})^{-1}\right)$ is the  Mille Cr\^epes coordinate chart  of $R_{s,p,n}^{-1}(U_0)$ with $\tau=\left(\begin{matrix}
	1&2&\cdots&r\\
	s+1&s+2&\cdots&s+r\\
\end{matrix}\right)\in\mathbb J_0$. 

Denote by $\iota:D^-_1\cap\mathcal A^{\tau}\hookrightarrow\mathcal A^{\tau}$ the natural embedding. Define a morphism $\mathcal P$ from $\mathcal A^{\tau}={\rm Spec}\,\mathbb Z[\overrightarrow A,\widetilde Y,\overrightarrow B^1,\cdots,\overrightarrow B^r]$ to 
$D^-_1\cap \mathcal A^{\tau}={\rm Spec}\,\mathbb Z[b_{2(s+2)},\cdots,b_{r(s+r)},\widetilde Y,\overrightarrow B^1,\cdots,\overrightarrow B^r]$
by the natural projection omitting $b_{1(s+1)}$. Computation yields that the following diagram commutes.
\vspace{-0.05in}
\begin{equation*}
\begin{tikzcd}
&\mathcal A^{\tau}\arrow[d,bend right=15,"{\mathcal P}"']\arrow{rr}{\mathfrak P}\,\,&&\mathbb {P}^{N^0_{s,p,n}}\times\cdots\times\mathbb {P}^{N^r_{s,p,n}}\\
&D^-_1\cap\mathcal A^{\tau}\arrow[u,bend right=15,"{\iota}"']\arrow{urr}{\mathfrak P}&& \\
\end{tikzcd}\vspace{-20pt}\,.
\end{equation*}
Then Lemma \ref{iemb} follows.\,\,\,\,$\endpf$.
\medskip

Denote by $\check {\mathfrak P}:D^-_1\rightarrow \mathfrak P\left(D^-_1\right)=\mathcal M_{s,p,n}$ the induced isomorphism and by $(\check {\mathfrak P})^{-1}$ its inverse. Define a map  $\mathcal P_{s,p,n}:\mathcal T_{s,p,n}\rightarrow D^-_1$ by
\begin{equation*}
\mathcal P_{s,p,n}:=(\check {\mathfrak P})^{-1}\circ \mathfrak P.
\end{equation*}
Immediately we have 
\begin{corollary}\label{retra}
$\mathcal P_{s,p,n}:\mathcal T_{s,p,n}\rightarrow D^-_1$ is a retraction, that is, the restriction of $\mathcal P_{s,p,n}$ to $D^-_1$ is the identity. 
\end{corollary}

\begin{lemma}\label{flat}
$\mathcal P_{s,p,n}:\mathcal T_{s,p,n}\rightarrow D^-_1$ is flat. 
\end{lemma}
{\bf \noindent Proof of Lemma \ref{flat}.}
It suffices to prove that $\mathfrak P:\mathcal T_{s,p,n}\rightarrow\mathcal M_{s,p,n}$ is flat. 

Take an arbitrary Mille Cr\^epes coordinate chart $\left(\mathcal A^{\tau},(J_l^{\tau})^{-1}\right)$ of $R_{s,p,n}^{-1}(U_l)$, $0\leq l\leq r$, with $\tau=\left(\begin{matrix} i_1&i_2&\cdots&i_r\\ j_1&j_2&\cdots&j_r\\
\end{matrix}\right)\in\mathbb J_l$.
Computation yields that when $l=0$, $\mathfrak P|_{\mathcal A^{\tau}}$ is equivalent to the product of the identity morphism of $\mathbb A^{n(n-p)-1}$ and the structure morphism ${\rm Spec}\,\mathbb Z[b_{i_1j_1}]\rightarrow{\rm Spec}\,\mathbb Z$; when $l=r$, $\mathfrak P|_{\mathcal A^{\tau}}$ is equivalent to the product of the identity morphism of $\mathbb A^{n(n-p)-1}$ and the structure morphism ${\rm Spec}\,\mathbb Z[a_{i_1j_1}]\rightarrow{\rm Spec}\,\mathbb Z$; when $1\leq l\leq r-1$, $\mathfrak P|_{\mathcal A^{\tau}}$ is equivalent to the product of the identity morphism of $\mathbb A^{n(n-p)-2}$ and the morphism ${\rm Spec}\,\mathbb Z[b_{i_1j_1},a_{i_{r-l+1}j_{r-l+1}}]\rightarrow{\rm Spec}\,\mathbb Z[t]$ defined by $t\mapsto  a_{i_{r-l+1}j_{r-l+1}}b_{i_1j_1}.$

We complete the proof of Lemma \ref{flat}.\,\,\,$\endpf$



\begin{lemma}\label{sti}
$\mathcal P_{s,p,n}$ induces an isomorphsim from $D_1^+$ to $D_1^-$.
\end{lemma}

{\noindent\bf Proof of  Lemma \ref{sti}.} 
The proof is the same as that for Lemma \ref{iemb}. 
\begin{remark}\label{gfi}
There is a geometric interpretation of  $\mathcal P_{s,p,n}$ as follows. According to Remark \ref{nd}, for each point on $\mathcal T_{s,p,n}$ we can flow it back to the source 
through  a sequence of adjacent orbital curves.  The isomorphism between $D^-_1$ and $D^+_1$ in Theorem \ref{gwond} is given by such flow.
\end{remark}

Recall the following realization of 
the spaces of complete collineations (see \cite{Pez} for instance).
\begin{example}\label{vcc}
Let $\mathbb P(M_{p\times q})$ be the projectivization of $M_{p\times q}$ the set of all $p\times q$ matrices.
It is clear that $\mathbb P(M_{p\times q})$ is ${\rm GL}_p\times {\rm GL}_q$-invariant under the left multiplication by ${\rm GL}_p$ and the right multiplication by ${\rm GL}_q$. Define $r:=\min\left\{p,q\right\}$.
Then $\mathbb P(M_{p\times q})$ has $r$ ${\rm GL}_p\times {\rm GL}_q$-orbits, whose closures $Z_{r}\supset Z_{r-1}\supset \cdots\supset Z_1$ are given by the condition
that $Z_i$ is the set of points corresponding to matrices of rank at most $i$.  Blowing up $\mathbb P(M_{p\times q})$ successively along the strict transform of
$Z_1,\cdots,Z_{r-1}$, we obtain the spaces of complete collineations $\widetilde {\mathbb 
 P}(M_{p\times q})$.      
\end{example}

{\bf\noindent Proof of Theorem \ref{red2}.} By the  Mille Cr\^epes coordinates introduced in \S \ref{vanderp},  we can show that $\widetilde {\mathbb P}(M_{p\times q})$ is isomorphic to $D^+_1$. By Lemmas \ref{emb}, \ref{iemb}, \ref{sti}, we have $D^+_1\cong\mathcal M_{q,p,p+q}$. The proof is complete. \,\,\,$\endpf$

\medskip

{\noindent\bf Proof of  Theorem \ref{gwond}.} 
By Lemma \ref{cpc}, $R_{s,p,n}:\mathcal T_{s,p,n}\rightarrow G(p,n)$ extends $\mathcal K^{-1}_{s,p,n}$. By Proposition \ref{tspnsmooth}, $\mathcal T_{s,p,n}$ is smooth. 

By (\ref{tschubert}) and Lemma \ref{snc}, we conclude that the complement of the open
${\rm GL}_s\times {\rm GL}_{n-s}$-orbit in $\mathcal T_{s,p,n}$ is a simple normal crossing divisor which  consists of $2r$ smooth, irreducible divisors 
$D^-_1, D^-_2,\cdots,D^-_r,D^+_1,D^+_2,\cdots,D^+_r$. Let $\left(\mathcal A^{\tau},(J^{\tau}_l)^{-1}\right)$ be any Mille Cr\^epes coordinate chart of $R_{s,p,n}^{-1}(U_l)$ with $0\leq l\leq r$ and $\tau=\left(\begin{matrix} i_1&i_2&\cdots&i_r\\ j_1&j_2&\cdots&j_r\\
\end{matrix}\right)\in\mathbb J_l$. 
Then $\mathfrak a, \mathfrak b\in \mathcal A^{\tau}$ are in the same orbit if and only if for $1\leq m\leq r-l$, \begin{equation*}
b_{i_mj_m}\in\mathfrak b\Longleftrightarrow b_{i_mj_m}\in\mathfrak a,\\
\end{equation*}
and for $r-l+1\leq m\leq r$,
\begin{equation*}
a_{i_mj_m}\in\mathfrak b\Longleftrightarrow a_{i_mj_m}\in\mathfrak a.\\
\end{equation*}
Hence the
${\rm GL}_s\times {\rm GL}_{n-s}$-orbit of $\mathcal T_{s,p,n}$ one to one corresponds to  (\ref{inrule}). It is clear that the closure of each orbit is smooth. Therefore, Property (C) follows.


Property (A) follows from Lemmas \ref{emb}, \ref{iemb}, \ref{sti}.

It is clear that the  complement of the open ${\rm GL}_s\times {\rm GL}_{n-s}$-orbit of $D^-_1\cong\mathcal M_{s,p,n}$ is the union of 
$\check  D_2,\,\check D_3,\, \cdots,\, \check D_r$. As in the case of the spaces of complete collineations (see Example \ref{vcc} for instance), we can show that $D^-_1\cong\mathcal M_{s,p,n}$ is wonderful by exploiting the Mille Cr\^epes coordinate charts as above. Hence Property (D) holds.

By Lemma \ref{flat} $\mathcal P_{s,p,n}:\mathcal T_{s,p,n}\rightarrow D_1^-$  is ${\rm GL}_s\times {\rm GL}_{n-s}$-equivariant and flat. $\mathcal P_{s,p,n}$ is a retraction by Corollary \ref{retra}, and the restriction of $\mathcal P_{s,p,n}$ on $D^+_1$ is an isomorphism by Lemma \ref{sti}. 
Noticing that the $\mathbb G_m$-action on $\mathcal M_{s,p,n}$ is the identity, we conclude  by (\ref{11}), (\ref{tschubert}) that $\mathcal P_{s,p,n}(\mathcal D_{(p-r-1+i,r+1-i)})=\mathcal P_{s,p,n}(D^+_i)=\mathcal P_{s,p,n}(D^-_{r+2-i})$ for $2\leq i\leq r$.


We complete the proof of Theorem \ref{gwond}. \,\,\,$\endpf$
\medskip

To prove Theorem \ref{moduli}, we first prove the following lemma.
\begin{lemma}\label{moduli1}
Base change (\ref{KF}) to an algebraically closed field $\mathbb K$. The following holds. 
\begin{enumerate}[label=(\alph*)]

\item For each closed point $q\in \mathcal M_{s,p,n}$, the fiber $Z_q:=(\mathcal P_{s,p,n})^{-1}(q)$ consists of a chain of $\mathbb G_m$-stable smooth rational curves.

\item The restriction of $R_{s,p,n}$ to $Z_q$ is an embedding.

\item Generically $R_{s,p,n}(Z_q)$ is a smooth rational curve of degree $r$ with respect to the Pl\"ucker embedding.

\item For any closed points $q,q^{\prime}\in \mathcal M_{s,p,n}$,  $R_{s,p,n}(Z_q)=R_{s,p,n}(Z_{q^{\prime}})$ if and only if $q=q^{\prime}$.


\end{enumerate} 
\end{lemma}
{\noindent\bf Proof of Lemma \ref{moduli1}.} 
Take a Mille Cr\^epes coordinate chart $\mathcal A^{\tau}$ of $R^{-1}_{s,p,n}(U_0)$. It is clear that the fiber $Z_q$ of $\mathcal P_{s,p,n}$ over a generic point $q$ is defined in $\mathcal A^{\tau}$ by fixing all the variables except $a_{i_1j_1}$. Note that the image of $Z_q\cap\mathcal A^{\tau}$ is given by the projection to $\mathbb P^{N_{p,n}}$. Plugging into (\ref{f1}), (\ref{f2}), (\ref{f3}),  each component is a polynomial in $a_{i_1j_1}$ such that the highest degree of $a_{i_1j_1}$ in the components are $r$ and that there is one component taking value $1$. Property (c) follows. 

Take any Mille Cr\^epes coordinate chart $\mathcal A^{\tau}$ of $R^{-1}_{s,p,n}(U_l)$ with $0\leq l\leq r$. Computation yields that any fiber $Z_q$ restricted to $\mathcal A^{\tau}$ is defined by fixing $a_{i_1j_1}\cdot b_{i_{r-l+1}j_{r-l+1}}$ and all the other variables. Then it is easy to verify that $Z_q$ is $1$-dimensional, reduced, and consists of a chain of $\mathbb G_m$-stable smooth rational curves, and, moreover, the restriction of $\phi$ to $Z_q$ is an embedding.
We conclude Properties (a), (b).

Next suppose that $R_{s,p,n}(Z_q)=R_{s,p,n}(Z_{q^{\prime}})$. We can write $Z_q$ as a chain of a rational curves $\cup_{k=1}^m\gamma_k$ with integers $1=k_0\leq k_1<\cdots<k_m=r$ such that 
$$\gamma_{1}\subset \overline{\mathcal D_{(p-k_1,k_1)}^-},\,\,\gamma_{i}\subset \overline{\mathcal D_{(p-k_i,k_i)}^-}\cap\overline{\mathcal D_{(p-k_{i-1},k_{i-1})}^+}\,\,{\rm for\,\,}2\leq i\leq m-1,\,\,{\rm and\,\,}\gamma_{m}\subset \overline{\mathcal D_{(p-k_{m-1},k_{m-1})}^+}.$$
Notice that for any fixed $i=1,\cdots,m$, the rational map $f_s^k$ given by (\ref{fsk}) is well-defined on $\gamma_i$ for  $k_{i-1}\leq k\leq k_i$. We can thus conclude that $q=q^{\prime}$. Property (d) follows. 

We complete the proof of Lemma \ref{moduli1}. \,\,\,$\endpf$
\medskip

{\noindent\bf Proof of Theorem \ref{moduli}.} By Property (B) in Theorem \ref{gwond} and Properties (b), (d) in Lemma \ref{moduli1}, there is a bijective morphism $\mathcal M_{s,p,n}\longrightarrow G(p,n)/ \! \! /\mathbb G_m$. 
We next show that $G(p,n)/ \! \! /\mathbb G_m$ is smooth following the idea in the proof of \cite[Proposition 6.3]{ORCW}.

Denote by $C_g$ the closure of the general orbit of the action in $G(p,n)$, and by $[C_g]$ the corresponding point in the Hilbert scheme ${\rm Hilb}(G(p,n))$. Notice that the tangent bundle $TG(p,n)$ of $G(p,n)$ is globally generated. Then  $H^1(\mathbb P^1,\mu^*TG(p,n))=0$ for any morphism $\mu:\mathbb P^1\to G(p,n)$, and 
$H^1(C_g,\mathcal N_{C_g/G(p,n)})=H^1(C_g,{\left.T_{G(p,n)}\right|}_{C_g})=0$ by Property (c) in Lemma \ref{moduli}.  Hence ${\rm Hilb}(G(p,n))$ is smooth at $[C_{g}]$, and there is a unique irreducible component of ${\rm Hilb}(G(p,n))$ containing it, which we denote by $\mathcal H$. Consider the induced $\mathbb G_m$-action on $\mathcal H$. By \cite[Theorem 2.5]{Bi}, there is a unique connected component of the fixed point scheme $\mathcal M \subset \mathcal H$ containing $[C_g]$ such that $\mathcal M$ is smooth at $[C_g]$. It is clear that the generic point of $\mathcal M$ parametrizing a smooth $\mathbb G_m$-invariant rational curve passing through the source $\mathcal V_{p,0}$ and the sink $\mathcal V_{p-r,r}$. Then $\mathcal M\cong G(p,n)/ \! \! /\mathbb G_m$. Now by \cite[Theorem 2.5]{Bi} again,
it suffices to show that $\mathcal H$ is smooth at every point of $\mathcal M$.

Since the morphism $\mathcal M_{s,p,n}\longrightarrow G(p,n)/ \! \! /\mathbb G_m$ is surjective, for every $[C]\in \mathcal M$, the corresponding subscheme $C$ is a chain of smooth rational curves with simple normal crossings by Property (a) in Lemma \ref{moduli1}. Hence the conormal sheaf $\mathcal I_{C,G(p,n)}/\mathcal I^2_{C,G(p,n)}$ is locally free. Applying the functor ${\rm Hom}(\,\cdot\,,\mathcal O_C)$ to the short exact sequence
$$
0\rightarrow{\mathcal I_{C,G(p,n)}/\mathcal I^2_{C,G(p,n)}}\rightarrow{\left.\Omega_{G(p,n)}\right|}_{C}\rightarrow{\Omega_C}\rightarrow0,$$
we get that ${\rm Ext}^1(\mathcal I_{C,G(p,n)}/\mathcal I^2_{C,G(p,n)},\mathcal O_C)$ is a quotient of ${\rm Ext}^1({\left.\Omega_{G(p,n)}\right|}_{C},\mathcal O_C)=H^1(C,{\left.T_{G(p,n)}\right|}_{C})$. Moreover, we can prove by induction that $H^1(C,{\left.T_{G(p,n)}\right|}_{C})=0$
as in the proof of \cite[Lemma~10]{FP}. Hence $\mathcal M$ is smooth at $[C]$ by \cite[Theorem~I.2.8]{Ko}.

We complete the proof of Theorem \ref{moduli}.\,\,\,$\endpf$


\section{Birational Morphisms among Rational Schemes}

\subsection{Kausz compactification and Landsberg-Manivel birational map}\label{kcpt}
We first recall the construction of Kausz \cite{Ka} with a slight generalization for the non-group compactification case.   We assume that $2p\leq n$ and take $s=n-p$ (hence $r=p$). 

Define $X^{(0)}:={\rm Proj}\,\mathbb Z[x_{00}, x_{ij} (1\leq i\leq p,1\leq j\leq n-p)]$, and its closed subschemes
\vspace{-.08in}
\begin{equation}\label{lblow}
\begin{tikzcd}
&Y^{(0)}_0\ar[hook]{r}&Y^{(0)}_1\ar[hook]{r}&\cdots\ar{r}&Y^{(0)}_{p-1}&\\
&&Z^{(0)}_{p-1}\ar[hook]{u}\ar[hook]{r}&\cdots\ar[hook]{r}&Z^{(0)}_1\ar[hook]{u}\ar[hook]{r}&Z^{(0)}_0\\
\end{tikzcd}.\vspace{-20pt}
\end{equation}
by $Y^{(0)}_l:=V(\mathcal I^{(0)}_l)$, $Z^{(0)}_l:=V(\mathcal J^{(0)}_l)$, $0\leq l\leq p-1$. Here $\mathcal I^{(0)}_l$ is the homogeneous ideal in $\mathbb Z[x_{00}, x_{ij} (1\leq i\leq p,1\leq j\leq n-p)]$ generated by all $(l+1)\times(l+1)$-subdeterminants of the matrix $(x_{ij})_{1\leq i\leq p,1\leq j\leq n-p}$ for $0\leq l\leq p-1$, and  $\mathcal J^{(0)}_0=(x_{00})$, $\mathcal J^{(0)}_l:=(x_{00})+\mathcal I^{(0)}_{p-l}$ for $1\leq l\leq p-1$. For $1\leq k \leq p$, inductively define 
$X^{(k)}\rightarrow X^{(k-1)}$ to be the blow up of $X^{(k-1)}$ along the closed subscheme $Y^{(k-1)}_{k-1}\cup Z^{(k-1)}_{p-k}$. Define the subscheme $Y^{(k)}_{k-1}\subset X^{(k)}$ (resp.~$Z^{(k)}_{p-k}\subset X^{(k)}$) to be the inverse image of $Y^{(k-1)}_{k-1}$ (resp.~$Z^{(k-1)}_{p-k}$) under the morphism $X^{(k)}\rightarrow X^{(k-1)}$, and for $l\neq k-1$ (resp.~$l\neq p-k$) the subscheme $Y^{(k)}_l\subset X^{(k)}$ (resp.~$Y^{(k)}_l\subset X^{(k)}$) to be the strict transform of $Y^{(k-1)}_l\subset X^{(k-1)}$ (resp. $Z^{(k-1)}_l\subset X^{(k-1)}$). Set ${\rm KA}_{p,n}:=X^{(p)}$, $Y_l:=Y^{(p)}_l$, $Z_l:=Z^{(p)}_l$, and denote by
\begin{equation}\label{kmorp}
\mathcal {KA}:{\rm KA}_{p,n}\longrightarrow X^{(0)}\cong\mathbb P^{p(n-p)}
\end{equation} the corresponding blow-up. Note that when $n=2p$, ${\rm KA}_{p,n}$ is Kausz's modular compactification of ${\rm GL}_p$. 

As in \cite[\S 4]{Ka}, we can explicitly parametrize ${\rm KA}_{p,n}$ by affine schemes as follows. Let 
$(\alpha,\beta)\in S_p\times S_{n-p}$
and set
$$
x_{ij}^{(0)}(\alpha,\beta):= \frac{x_{\alpha(i),\beta(j)}}{x_{00}}
\qquad
(1\leq i\leq p,1\leq j\leq n-p).
$$
For $1\leq k\leq p$ we define elements
\begin{eqnarray*}
&y_{ji}(\alpha,\beta)\ (1\leq i\leq k,\, i<j\leq p),\,\,\,\, 
z_{ij}(\alpha,\beta)\
(1\leq i\leq k,\,  i<j\leq n-p),\\
&x_{ij}^{(k)}(\alpha,\beta)\ 
(k+1\leq i\leq p,\,  k+1\leq j\leq n-p) 
\end{eqnarray*}
of the function field 
$\mathbb Q(X^{(0)})=\mathbb Q(x_{ij}/x_{00}\ (1\leq i\leq p,\,1\leq j\leq n-p))$
inductively by
\begin{eqnarray*}
&y_{ik}(\alpha,\beta):= 
\frac{x_{ik}^{(k-1)}(\alpha,\beta)}
{x_{kk}^{(k-1)}(\alpha,\beta)} 
\ \ (k+1\leq i\leq p),\,\,\,\,
z_{kj}(\alpha,\beta):= 
\frac{x_{kj}^{(k-1)}(\alpha,\beta)}
{x_{kk}^{(k-1)}(\alpha,\beta)} \ \  (k+1\leq j\leq n-p),
\\
&x_{ij}^{(k)}(\alpha,\beta):= 
\frac{x_{ij}^{(k-1)}(\alpha,\beta)}
{x_{kk}^{(k-1)}(\alpha,\beta)}
-
y_{ik}(\alpha,\beta)\ 
z_{kj}(\alpha,\beta)\ \ \
(k+1\leq i\leq p,\,  k+1\leq j\leq n-p).
\end{eqnarray*}
Finally, we set $t_0(\alpha,\beta) := t_0 := x_{00}$ and
$$
t_i(\alpha,\beta) := 
t_0 \cdot
\prod\nolimits_{j=1}^i
x_{jj}^{(j-1)}(\alpha,\beta)
\ \ 
(1\leq i\leq p).\ \ \ \ 
$$
Define for each $l\in\{0,\dots,n\}$ a bijection
$\iota_l:\{1,\dots,n+1\} \xmapsto{\,\,} \{0,\dots,n\}$  by 
$$
\iota_l(i)=
\left\{
\begin{array}{lll}
i & \text{if\quad $1\leq i\leq l$} \\
0 & \text{if\quad $i=l+1$} \\
i-1 & \text{if\quad $l+2\leq i\leq n+1$}
\end{array}
\right.,\,\,\,{\rm for\,\,}1\leq i\leq n+1.
$$
With this notation, we define for each triple 
$(\alpha,\beta,l)\in S_n\times S_n\times [0,n]$
polynomial subalgebras $R(\alpha,\beta,l)$ of $\mathbb Q({\rm {KA}}_{p,n})=\mathbb Q(X^{(0)})$ by
\begin{equation}\label{kcoor}
\mathbb Z\left[
\frac{t_{\iota_l(i+1)}(\alpha,\beta)}{t_{\iota_l(i)}(\alpha,\beta)} \
(1\leq i\leq n),\
y_{ji}(\alpha,\beta)(1\leq i<j\leq p),\
z_{ij}(\alpha,\beta) \
(1\leq i\leq p,\,i<j\leq n-p)
\right].    
\end{equation}
Then as in \cite[Proposition 4.1]{Ka}, one can show that there is a covering of ${\rm KA}_{p,n}$ by open affine pieces
$X(\alpha,\beta,l)$\ 
$((\alpha,\beta,l)\in S_n\times S_n\times [0,n])$,
such that 
$\Gamma(X(\alpha,\beta,l),\mathcal O)=R(\alpha,\beta,l)
$ as subrings of the function field $\mathbb Q({\rm KA}_{p,n})$. Denote the corresponding morphism by
\begin{equation}\label{habl}
h_{(\alpha,\beta,l)}:{\rm Spec\,}R(\alpha,\beta,l)
\longrightarrow{\rm KA}_{p,n}.   
\end{equation}

We next recall the Landsberg-Manivel birational map \cite{LM} for Grassmannians in terms of the Pl\"ucker coordinate functions. Define a rational map $LM$ from  $X^{(0)}\cong\mathbb P^{p(n-p)}$ to the open subscheme $U_p\subset G(p,n)$ (see (\ref{up1})) by
\begin{equation*}    
[x_{00}, x_{ij} (1\leq i\leq p,1\leq j\leq n-p)]\xdashmapsto{\,\,\,}\left(\begin{matrix}
x_{11}&x_{12}&\cdots&x_{1(n-p)}\\
x_{21}&x_{22}&\cdots&x_{2(n-p)}\\
\cdots&\cdots&\cdots&\cdots\\
x_{p1}&x_{p2}&\cdots&x_{p(n-p)}
\end{matrix}\hspace{-.12in} \begin{matrix} &\hfill\tikzmark{c}\\ 
&\\
&\\
&\hfill\tikzmark{d} \end{matrix} \hspace{.12in} \begin{matrix}
 x_{00}&0&\cdots&0\\
  0&x_{00}&\cdots&0\\
\cdots&\cdots&\cdots&\cdots\\
0&0&\cdots&x_{00}   
\end{matrix}\right).
\tikz[remember picture,overlay]   \draw[dashed,dash pattern={on 4pt off 2pt}] ([xshift=0.5\tabcolsep,yshift=6pt]c.north) -- ([xshift=0.5\tabcolsep,yshift=3pt]d.south);
\end{equation*}
Up to a certain permutation of the homogeneous coordinates, the Landsberg-Manivel birational map $\mathcal {LM}:\mathbb P^{p(n-p)}\xdashrightarrow{\,\,\,}
G(p,n)$ is given by 
$$e\circ LM=[\cdots,P_I\circ LM,\cdots]_{I\in\mathbb I_{p,n}},$$
where $e:G(p,n)\xhookrightarrow{\,\,\,\,\,\,\,\,\,\,\,\,\,\,}\mathbb P^{N_{p,n}}$ is the Pl\"ucker embedding. 

Applying the Landsberg-Manivel birational map, we prove that
\begin{proposition}\label{g2kausz}
${\rm KA}_{p,n}\cong\mathcal T_{n-p,p,n}$ over ${\rm Spec}\,\mathbb Z$. In particular, Theorem \ref{K=T} holds.
\end{proposition}
{\noindent\bf Proof of Proposition \ref{g2kausz}.}
Computing in the coordinates (\ref{kcoor}), it is clear  that the rational $\mathcal {LM}\circ \mathcal {KA}:{\rm KA}_{p,n}\dashrightarrow G(p,n)$ extends to a regular morphism. Now consider 
\begin{equation}
\mathcal K_{n-p,p,n}\circ\mathcal {LM}\circ \mathcal {KA}:{\rm KA}_{p,n}\xdashrightarrow{\,\,\,}\mathbb {P}^{N_{p,n}}\times\mathbb {P}^{N^0_{s,p,n}}\times\cdots\times\mathbb {P}^{N^p_{s,p,n}}.     
\end{equation}
Similar to the proof of Lemma \ref{em}, we can show that $\mathcal K_{n-p,p,n}\circ\mathcal {LM}\circ \mathcal {KA}$ extends to an embedding. Therefore $\mathcal K_{n-p,p,n}\circ\mathcal {LM}\circ \mathcal {KA}$ extends to an isomorphism from ${\rm KA}_{p,n}$ to $\mathcal T_{n-p,p,n}$.

Proposition \ref{g2kausz} follows.\,\,\,$\endpf$
\medskip

From the proof of Proposition \ref{g2kausz}, it immediately follows that
\begin{proposition}\label{g3kausz}
We have 
the following commutative diagram.
\vspace{-0.05in}
\begin{equation}\label{K=F3}
\begin{tikzcd}
&&\,\,\mathcal T_{n-p,p,n}\arrow{dl}[swap]{\hspace{-.03in}\mathcal {KA}}\arrow{dr}{\hspace{-0.03in}R_{n-p,p,n}}\,\,&\\
&\mathbb P^{p(n-p)}\arrow[dashed]{rr}{\mathcal {LM}}&& G(p,n) \\
\end{tikzcd}\vspace{-20pt}\,.
\end{equation}
In particular, Corollary \ref{2kausz} holds.
\end{proposition}
\begin{remark}\label{KGE2}
$Y_i$ corresponds to $D^-_{i+1}$ and $Z_j$ corresponds to $D^+_{j+1}$ as indicated by (\ref{KGE}), (\ref{lblow}).
\end{remark}
\begin{remark}
A nice observation of \cite{Di} is that it suffices to blow up $Z_j$ and blow down $D^+_i$. It is also clear from our picture that such half data suffices, for $Y_i$ and $D^-_{i+1}$ are the same determinantal variety when restricted to the common affine subset $\mathbb A^{p(n-p)}$ of the projective space and of the Grassmannian.
\end{remark}

\subsection{Fibrations}
\label{basicif}

Next we turn to the fibration structures on $\mathcal T_{s,p,n}$ and $\mathcal M_{s,p,n}$. \medskip

{\bf\noindent Proof of Proposition \ref{tfiber}.} 
Denote by 
\begin{equation*}
\begin{split}
&p_0:\mathbb {P}^{N_{p,n}}\times\mathbb {P}^{N^0_{s,p,n}}\times\cdots\times\mathbb {P}^{N^r_{s,p,n}}\longrightarrow\mathbb {P}^{N^0_{s,p,n}},\,\,\,\,p_r:\mathbb {P}^{N_{p,n}}\times\mathbb {P}^{N^0_{s,p,n}}\times\cdots\times\mathbb {P}^{N^r_{s,p,n}}\longrightarrow\mathbb {P}^{N^r_{s,p,n}}\\
\end{split}   
\end{equation*}
the natural projections. It is clear that when $p<n-s$ and $p<s$, the restriction of $p_r$ to $\mathcal T_{s,p,n}$ is a morphism to $G(p,n-s)$. Fix any $0\leq l\leq r$ and $\tau=\left(\begin{matrix}
i_1&\cdots&i_r\\
j_1&\cdots&j_r\\
\end{matrix}\right)\in\mathbb J_l$. Take the open subscheme $U_{(1,2,\cdots,l,\delta_1,\delta_2,\cdots,\delta_{r-l})}$ of $G(p,n-s)$ as in (\ref{u12}) with $$\{\delta_1,\delta_2,\cdots,\delta_{r-l}\}=\{j_1-s,j_2-s,\cdots,j_{r-l}-s\}.$$  Let  $U_l\subset G(p,n)$, $\Gamma_l^{\tau}:\mathbb A^{p(n-p)}\rightarrow U_l$  be defined by (\ref{ul}), (\ref{ws}), respectively. Define a new index $\widetilde\tau=\left(\begin{matrix}
i_1&\cdots&i_{r-l}&i_{r-l+1}&\cdots&i_r\\
j^{\prime}_1&\cdots&j^{\prime}_{r-l}&j_{r-l+1}&\cdots&j_r\\
\end{matrix}\right)$ such that $(j^{\prime}_1,\cdots,j^{\prime}_{r-l})$ is a permutation of $(l+1,\cdots,r)$ and that for any $1\leq\alpha<\beta\leq r-l$, $j^{\prime}_{\alpha}<j^{\prime}_{\beta}$ if and only if $j_{\alpha}<j_{\beta}$. 
We now define an embedding from $U_{(1,2,\cdots,l,\delta_1,\delta_2,\cdots,\delta_{r-l})}\times\mathcal A^{\widetilde\tau}$ to $\mathcal A^{\tau}$ as follows.  Define $\widetilde U_l\subset G(p,s+p)$, $\widetilde \Gamma_l^{\tau}:\mathbb A^{ps}\rightarrow U_l$ by setting $n=s+p$ in (\ref{ul}), (\ref{ws}), respectively. Write matrices $U_l=(A\,\,B)$, $\widetilde U_l=(C\,\,D)$, where $A$, $B$, $C$, $D$ are $p\times s$, $p\times (n-s)$, $p\times s$, $p\times p$ matrices, respectively.  Define a morphism $\eta:U_{(1,2,\cdots,l,\delta_1,\delta_2,\cdots,\delta_{r-l})}\times\widetilde U_l\longrightarrow U_l$ by
\begin{equation*}
A=C,\,\, B=D\cdot U_{(1,2,\cdots,l,\delta_1,\delta_2,\cdots,\delta_{r-l})}.  
\end{equation*}
Computation yields that the birational map $(\Gamma_l^{\tau})^{-1}\circ({\rm Id},\widetilde\Gamma_l^{\tau}):\mathbb A^{p(n-s-p)}\times\mathbb A^{ps}\rightarrow U_l$ extends to an isomorphism.  We thus conclude that $\mathcal T_{s,p,n}$ is a locally trivial fibration over $G(p,n-s)$ with the fiber  $\mathcal T_{p,p,s+p}$.

The remaining cases can be proved in the same way. 
\,\,\,$\endpf$
\medskip




In the following, we produce explicit birational transformations among projective bundles over Grassmannians as follows. Let $E=E_1\oplus E_2$. Denote by $Gr(p,E_1)$ the sub-Grassmannian of $Gr(p,E)$
consisting of $p$-planes  in $E_1$, and denote by $\mathbb P(N_1)$ the projectivization of the normal bundle $N_1$ of $Gr(p,E_1)$ in $Gr(p,E)$. When $p\leq n-s$, define $Gr(p,E_2)$ in the same way; otherwise,  
define $Gr(p,E_2)$ to be the sub-Grassmannian of $Gr(p,E)$
consisting of $p$-planes $W$ such that $\dim W\cap E_1=p+s-n$. Denote by $\mathbb P(N_2)$  the projectivization of the normal bundle $N_2$ of $Gr(p,E_2)$ in $Gr(p,E)$. 

Restricting (\ref{sblow}) to $U_0$ with $\sigma$ given by $\sigma(k)=k+1$ for $0\leq k\leq p-1$ and  $\sigma(p)=0$, 
we derive the following commutative diagram 
\vspace{-.05in}
\begin{equation*}
\begin{tikzcd}
&D^-_1\arrow{d}{\hspace{-.6in} R_{s,p,n}\big|_{D^-_1}\,\,\,\,\,\,}\arrow{r}{\,\,\,\left(g^{\sigma}_1\circ\cdots\circ g^{\sigma}_{p}\right)\big|_{D^-_1}\,\,\,}&[4em]\mathbb P(N_1) \arrow{ld}{g^{\sigma}_{0}|_{\mathbb P(N_1)}} \\
&\hspace{-.7in}Gr(p,E_1)&\\
\end{tikzcd}\vspace{-20pt}\,.
\end{equation*}
Similarly, restricting (\ref{sblow}) to $U_r$ with a certain permutation $\sigma^{\prime}$, we have the following commutative diagram.
\vspace{-.05in}
\begin{equation*}
\begin{tikzcd}
&D^+_1\arrow{d}{\hspace{-.6in} R_{s,p,n}\big|_{D^+_1}\,\,\,\,\,\,}\arrow{r}{\,\,\,\left(g^{\sigma^{\prime}}_1\circ\cdots\circ g^{\sigma^{\prime}}_{p}\right)\big|_{D^+_1}\,\,\,}&[4em]\mathbb P(N_2) \arrow{ld}{g^{\sigma^{\prime}}_{0}|_{\mathbb P(N_2)}} \\
&\hspace{-.7in}Gr(p,E_2)&\\
\end{tikzcd}\vspace{-20pt}\,.
\end{equation*}
Denote by $R_1:D^-_1\rightarrow\mathbb  P(N_1)$ and $R_2:D^+_1\rightarrow\mathbb  P(N_2)$ the corresponding blow-ups. Then based on the observation that the source and the sink of $\mathcal T_{s,p,n}$ are isomorphic,  we have 
\medskip




{\noindent\bf Proof of Corollary \ref{pbun}.} It follows from Lemma \ref{sti}.\,\,\,$\endpf$

\bibliographystyle{plain} 
\bibliography{ref}

\begin{thebibliography}{10}

\bibitem{Al}
A.~Alguneid.
\newblock Degeneration of space collineations.
\newblock {\em Proc. Egyptian Acad. Sci.}, 7:1--17 (1952), 1951.

\bibitem{Bi}
A.~Bia{\l}ynicki-Birula.
\newblock Some theorems on actions of algebraic groups.
\newblock {\em Ann. of Math. (2)}, 98:480--497, 1973.

\bibitem{BS}
A.~Bia{\l}ynicki-Birula and A.~Sommese.
\newblock Quotients by {${\bf C}^{\ast} $} and {${\rm SL}(2,{\bf C})$} actions.
\newblock {\em Trans. Amer. Math. Soc.}, 279(2):773--800, 1983.

\bibitem{BP}
P.~Bravi and G.~Pezzini.
\newblock Primitive wonderful varieties.
\newblock {\em Math. Z.}, 282(3-4):1067--1096, 2016.

\bibitem{BM}
M.~Brion and S.~Kumar.
\newblock {\em Frobenius splitting methods in geometry and representationtheory}, volume 231 of {\em Progress in Mathematics}.
\newblock Birkh\"{a}user Boston, Inc., Boston, MA, 2005.

\bibitem{BV}
W.~Bruns and U.~Vetter.
\newblock {\em Determinantal rings}, volume 1327 of {\em Lecture Notes in Mathematics}.
\newblock Springer-Verlag, Berlin, 1988.

\bibitem{DP}
C.~De~Concini and C.~Procesi.
\newblock Complete symmetric varieties.
\newblock In {\em Invariant theory ({M}ontecatini, 1982)}, volume 996, pages 1--44. Springer, Berlin, 1983.

\bibitem{DP2}
C.~De~Concini and C.~Procesi.
\newblock Complete symmetric varieties. {II}. {I}ntersection theory.
\newblock In {\em Adv. Stud. Pure Math.} Vol. 6, 481--513, North-Holland, Amsterdam, 1985.

\bibitem{Di}
C.~Ding.
\newblock Birational transformations on irreducible compact {H}ermitian symmetric spaces.
\newblock {\em Int. Math. Res. Not. IMRN}, (11):9266--9291, 2024.

\bibitem{Fa}
G.~Faltings.
\newblock Explicit resolution of local singularities of moduli-spaces.
\newblock {\em J. Reine Angew. Math.}, 483:183--196, 1997.

\bibitem{Fu}
A.~Fujiki.
\newblock Fixed points of the actions on compact {K}\"{a}hler manifolds.
\newblock {\em Publ. Res. Inst. Math. Sci.}, 15(3):797--826, 1979.

\bibitem{FP}
W.~Fulton and R.~Pandharipande.
\newblock Notes on stable maps and quantum cohomology.
\newblock volume~62 of {\em Proc. Sympos. Pure Math.}, pages 45--96. Amer. Math. Soc., Providence, RI, 1997.

\bibitem{He2}
X.~He.
\newblock Normality and {C}ohen-{M}acaulayness of local models of {S}himura varieties.
\newblock {\em Duke Math. J.}, 162(13):2509--2523, 2013.

\bibitem{Hw1}
J.~Hwang.
\newblock Nondeformability of the complex hyperquadric.
\newblock {\em Invent. Math.}, 120(2):317--338, 1995.

\bibitem{Kap}
M.~M. Kapranov.
\newblock Chow quotients of {G}rassmannians. {I}.
\newblock volume~16 of {\em Adv. Soviet Math.}, pages 29--110. Amer. Math. Soc., Providence, RI, 1993.

\bibitem{Ka}
I.~Kausz.
\newblock A modular compactification of the general linear group.
\newblock {\em Doc. Math.}, 5:553--594, 2000.

\bibitem{Kn}
F.~Knop.
\newblock The {L}una-{V}ust theory of spherical embeddings.
\newblock In {\em Proceedings of the {H}yderabad {C}onference on {A}lgebraic {G}roups}, pages 225--249. Manoj Prakashan, Madras, 1991.

\bibitem{Ko}
J\'{a}nos Koll\'{a}r.
\newblock {\em Rational curves on algebraic varieties}, volume~32 of {\em Ergebnisse der Mathematik und ihrer Grenzgebiete (3)}.
\newblock Springer-Verlag, Berlin, 1996.

\bibitem{La0}
D.~Laksov.
\newblock The geometry of complete linear maps.
\newblock {\em Ark. Mat.}, 26(2):231--263, 1988.

\bibitem{LLT}
D.~Laksov, A.~Lascoux, and A.~Thorup.
\newblock On {G}iambelli's theorem on complete correlations.
\newblock {\em Acta Math.}, 162(3-4):143--199, 1989.

\bibitem{LM}
J.~Landsberg and L.~Manivel.
\newblock Construction and classification of complex simple {L}ie algebras via projective geometry.
\newblock {\em Selecta Math. (N.S.)}, 8(1):137--159, 2002.

\bibitem{L}
D.~Luna.
\newblock Toute vari\'{e}t\'{e} magnifique est sph\'{e}rique.
\newblock {\em Transform. Groups}, 1(3):249--258, 1996.

\bibitem{Lu2}
G.~Lusztig.
\newblock Parabolic character sheaves. {I}.
\newblock {\em Mosc. Math. J.}, 4(1):153--179, 311, 2004.

\bibitem{ORCW}
G.~Occhetta, E.~Romano, L.~Conde, and J.~Wi\'sniewski.
\newblock Chow quotients of {$\mathbb C^*$}-action.
\newblock {\em arxiv:2310.18623}, 2023.

\bibitem{Pez}
G.~Pezzini.
\newblock Lectures on wonderful varieties.
\newblock {\em Acta Math. Sin. (Engl. Ser.)}, 34(3):417--438, 2018.

\bibitem{Se2}
J.~Semple.
\newblock The variety whose points represent complete collineations of {$S_r$} on {$S'_r$}.
\newblock {\em Univ. Roma Ist. Naz. Alta Mat. Rend. Mat. e Appl. (5)}, 10:201--208, 1951.

\bibitem{Sev2}
F.~Severi.
\newblock I fondamenti della geometria numerativa.
\newblock {\em Ann. Mat. Pura Appl. (4)}, 19:153--242, 1940.

\bibitem{Siu}
Y.~Siu.
\newblock Nondeformability of the complex projective space.
\newblock {\em J. Reine Angew. Math.}, 399:208--219, 1989.

\bibitem{Sp1}
T.~A. Springer.
\newblock Intersection cohomology of {$B\times B$}-orbit closures in group compactifications.
\newblock volume 258, pages 71--111. 2002.
\newblock With an appendix by Wilberd van der Kallen, Special issue in celebration of Claudio Procesi's 60th birthday.

\bibitem{Stu}
E.~Study.
\newblock Ueber die {G}eometrie der {K}egelschnitte, insbesondere deren {C}harakteristikenproblem.
\newblock {\em Math. Ann.}, 27(1):58--101, 1886.

\bibitem{Th1}
M.~Thaddeus.
\newblock Complete collineations revisited.
\newblock {\em Math. Ann.}, 315(3):469--495, 1999.

\bibitem{TK}
A.~Thorup and S.~Kleiman.
\newblock Complete bilinear forms.
\newblock In {\em Algebraic geometry ({S}undance, {UT}, 1986)}, volume 1311 of {\em Lecture Notes in Math.}, pages 253--320. Springer, Berlin, 1988.

\bibitem{Ty}
J.~A. Tyrrell.
\newblock Complete quadrics and collineations in {$S_n$}.
\newblock {\em Mathematika}, 3:69--79, 1956.

\bibitem{Va}
I.~Vainsencher.
\newblock Complete collineations and blowing up determinantal ideals.
\newblock {\em Math. Ann.}, 267(3):417--432, 1984.

\bibitem{Van}
B.~L. van~der Waerden.
\newblock Zur algebraischen {G}eometrie. {XV}.
\newblock {\em Math. Ann.}, 115(1):645--655, 1938.

\end{thebibliography}
\Addresses
\end{document}